\documentclass[a4paper,11pt,onecolumn,oneside,notitlepage,openbib,final]{article}
\pdfoutput=1
\usepackage{mathpazo}	
\usepackage{relsize}
\renewcommand{\ast}{{\mathlarger *}} 
\usepackage{float}
\usepackage{lscape}	
\usepackage[bottom]{footmisc}			
\usepackage[margin=1in]{geometry}
\usepackage{setspace} 
\usepackage{amsmath,amsfonts,amssymb,amsthm}
\usepackage{bbm} 
\usepackage{mathrsfs}		
\allowdisplaybreaks
\usepackage{thmtools}
\usepackage{mathtools}
\usepackage[pdflatex]{graphicx}
\usepackage{color}
\usepackage[labelformat=parens]{subfig}		
\usepackage[notcite,color]{showkeys}
\definecolor{refkey}{rgb}{0.6,0.4,1}
\definecolor{labelkey}{rgb}{0.6,0.4,1}
\usepackage[pdftex,bookmarks=true,bookmarksnumbered=true,hidelinks]{hyperref}
\usepackage{url}

\makeatletter
\let\orgdescriptionlabel\descriptionlabel
\renewcommand*{\descriptionlabel}[1]{%
	\let\orglabel\label
	\let\label\@gobble
	\phantomsection
	\edef\@currentlabel{#1\unskip}%
	\let\label\orglabel
	\orgdescriptionlabel{#1}%
}
\makeatother

\usepackage{natbib}
\usepackage{enumerate}
\usepackage{alltt}	
\usepackage{cleveref} 
\usepackage[most]{tcolorbox}
\usepackage{listings} 
\usepackage {spverbatim} 
\usepackage{adjustbox} 
\usepackage{multirow}
\usepackage{array}
\newcolumntype{C}{>{$\displaystyle} c <{$}}

\def\dif#1#2{\frac{d #1}{d #2}}

\def\pdif#1#2{\frac{\partial #1}{\partial #2}}

\def\tpdif#1#2{\partial #1 / \partial #2}

\def\2'{^{\prime\prime}}

\def\fB{\mathfrak B}
\def\fb{\mathfrak b}
\def\c{\mathbf c}
\def\d{\mathbf d}
\def\D{\mathbf D}
\def\E{\mathbb E}
\def\e{\mathbf e}
\def\f{\mathbf f}
\def\F{\mathbf F}
\def\g{\mathbf g}
\def\h{\mathbf h}

\def\N{\mathbb N}
\def\P{\mathbb P}
\def\p{\mathbf p}

\def\R{\mathbb R}

\def\V{\mathbf V}
\def\v{\mathbf v}
\def\w{\mathbf w}

\def\x{\mathbf x}
\def\y{\mathbf y}

\def\z{\mathbf z}
\def\0{\mathbf 0}
\def\1{\mathbf 1}
\def\O{\mathbf O}

\def\prA{A'}
\def\cA{\mathcal A}

\def\cX{\mathcal X}

\def\cZ{\mathcal Z}
\def\cP{\mathcal P}
\def\cQ{\mathcal Q}

\def\bpi{{\boldsymbol{\pi}}}
\def\btheta{{\boldsymbol{\theta}}}

\newcommand{\qandq}{\quad\text{and}\quad}

\DeclareMathOperator*{\argmax}{argmax}

\DeclareFontFamily{U}{mathx}{\hyphenchar\font45}
\DeclareFontShape{U}{mathx}{m}{n}{
      <5> <6> <7> <8> <9> <10>
      <10.95> <12> <14.4> <17.28> <20.74> <24.88>
      mathx10
      }{}
\DeclareSymbolFont{mathx}{U}{mathx}{m}{n}
\DeclareMathSymbol{\bigtimes}{1}{mathx}{"91}

\DeclareMathOperator{\NE}{NE}
\DeclareMathOperator{\RP}{RP}
\theoremstyle{plain}
\newtheorem{thm}{Theorem}\crefname{thm}{theorem}{theorems}
\newtheorem{lem}{Lemma}\crefname{lem}{lemma}{lemmas}
\newtheorem{prop}{Proposition}\crefname{prop}{proposition}{propositions}
\newtheorem{cor}{Corollary}\crefname{cor}{corollary}{corollaries}

\crefname{obs}{observation}{observations}

\theoremstyle{definition}
\newtheorem{dfn}{Definition}\crefname{dfn}{definition}{definitions}

\newtheorem*{assmp*}{Assumption}

\newtheorem{innerAssmp}{Assumption}\crefname{innerAssmp}{assumption}{assumptions}
\newenvironment{assmp}[1]
{\renewcommand\theinnerAssmp{#1}\innerAssmp}
{\endinnerAssmp}

\crefname{assmpF}{Assumption}{Assumptions}
\newtheorem{assmpQ}{Assumption}
\crefname{assmpQ}{Assumption}{Assumptions}

\newtheorem{assmpA}{Assumption}

\crefname{assmpA}{Assumption}{Assumptions}

\theoremstyle{remark}
\newtheorem*{rmk}{Remark}

\newtheorem{exmpl}{Example}\crefname{exmp}{example}{examples}

\crefname{figure}{figure}{figures}

\newcounter{exmp}
\newenvironment{exmp}[1][]{\stepcounter{exmp}\begin{exmpl}[#1]}{\hspace{-\baselineskip}\hfill \rule{0.33em}{0.8em}\end{exmpl}}

\usepackage[font=small,labelfont=bf]{caption}
\DeclareCaptionFormat{myMain}{#1#2#3\vspace{-0.5\baselineskip}\hrulefill}
\DeclareCaptionFormat{mySub}{#1#2#3}
\title{Net gains in evolutionary dynamics\\ \large A unifying and intuitive approach to dynamic stability\footnote{On top of any papers of mine, Bill Sandholm gave me numerous advice and constant encouragement for this paper over years. I heartily dedicate this paper to him. The paper was presented at Kyoto, Hitotsubashi, Temple, Tohoku, Tsukuba, Game Theory Society world congress, Stony Brook international game theory conference, Learning, Evolution and Games conference, Japan Economic Association meeting, Tohoku-UEA joint workshop, and Game Theory workshop. The author thanks to Takashi Akamatsu, Kajii Atsushi, Susumu Cato, Dimitrios Diamantaras, Shota Fujishima, Daisuke Oyama, Chiaki Hara, Daisuke Hirata, Ryo Itoh, Michi Kandori, Yusuke Kasuya, Tatsuhito Kohno, Jonathan Newton, Ryoji Sawa, Karl Schlag, Tadashi Sekiguchi, Jeff Shamma, Dao-Zhi Zeng, and the participants at these presentations for valuable comments. Declarations of interest: none.}
}
\author{Dai {\sc Zusai}\thanks{Graduate School of Economics and Management, Tohoku University, Sendai, Japan (current workplace); Department of Economics, Temple University, Philadelphia, U.S.A. (the place where most of the work was done.) E-mail: \texttt{ZusaiDPublic@gmail.com}.}}
\date{\today}

\begin{document}
\setcounter{page}{0}
\maketitle
\thispagestyle{empty}
\begin{abstract}

	Static stability in economic models means negative incentives for deviation from equilibrium strategies, which we expect to assure a return to equilibrium, i.e., dynamic stability, as long as agents respond to incentives. There have been many attempts to prove this link, especially in evolutionary game theory, yielding both negative and positive results. This paper presents a universal and intuitive approach to this link. We prove that static stability assures dynamic stability if agents' choices of switching strategies are rationalizable by introducing costs and constraints in those switching decisions. This idea guides us to define \textit{net }gains from switches as the payoff improvement after deducting the costs. Under rationalizable dynamics, an agent maximizes the expected net gain subject to the constraints. We prove that the aggregate maximized expected net gain works as a Lyapunov function.  It also explains reasons behind the known negative results. While our analysis here is confined to myopic evolutionary dynamics in population games, our approach is applicable to more complex situations. 
	
	
\end{abstract}
\newpage

\section{Introduction}\label{sec:intro}
Equilibrium stability has been a fundamental and classic issue in economic theory. While it needs a rigorous formulation of an off-equilibrium adjustment process to study dynamic stability, conventional economists typically believe in the idea of static stability: it basically means that deviation from an equilibrium generates an incentive to cancel the deviation. If this holds, economists conjecture that economic agents should return to the equilibrium. This idea is common to, for example, evolutionary stable state (ESS) in games and Hicksian stability in competitive market theory. Like the study of t\^{a}tonnement processes for the latter, evolutionary game theory has been devoted to establishing \textit{the link between static and dynamic stability} by considering a variety of evolutionary dynamics that reflect various forms of bounded rationality and off-equilibrium learning \citep{SandholmPopText,Young04}. From those studies, we do have both positive and negative results; see \citet{Friedman91Ecta_EvolGamesEcon} and \citet{HofSand09JET_StableGames}---henceforth H\&S. However, we have not found an ultimate \textit{economic principle} to complete the link and to distinguish negative results from positive results. The current lack of the ultimate principle prevents us from extending those results obtained in a simple setting of normal form games to  more complicated situations. Here in the current paper, we solve this remaining gap by presenting an economically intuitive approach to constructing a \textit{disequilibrium index }that tells whether a given dynamic converges or not.

For economists, it is natural to measure incentives for agents to switch strategies. Equilibrium is defined as a state where there is no such incentive, and we expect that agents settle down to equilibrium strategies as the incentive for a switch diminishes. Economic intuition guides us to calculate payoff improvement by a switch, i.e., the payoff difference between an old strategy and a new strategy, as the measure of incentive for a switch; we call it a \textit{gross gain}. However, it does not work generally as a predictor of evolutionary dynamics. Evolutionary dynamics are typically coined with bounded rationality. Agents in those dynamics may miss an opportunity to switch to a better strategy or may not take the best response strategy. So, gross gains are not a perfect driver to guide agents' strategy adjustments (\Cref{sec:GrossGain}).

This observation motivates us to finding the exact factor to tell when an agent chooses to switch the strategy and what strategy to switch to; this is what we do in this paper. We first reconstruct a variety of bounded rationality in various evolutionary dynamics with an explicit formulation of costs and constraints in revision of strategies (\Cref{sec:Basic}). That is, if an agent dismisses an opportunity to switch to a better strategy, we rationalize it by introducing a switching cost. If an agent switches to a suboptimal strategy, then we regard that there is a constraint on available strategies and the exact (unconstrained) optimal strategy was not available. Then, we define the true incentive for revision or a \textit{net gain} by calculating a payoff difference between the current strategy and the best strategy among available strategies and then subtracting a switching cost from it. The net gain explains exactly whether an agent switches a strategy \textit{or not} and which strategy the agent switches to, while the gross gain does not. To explain continuity in switching rates, we allow for stochastic switching costs  and probabilistic available strategy sets. We define an \textit{ex-ante} net gain by taking the expected value of net gains over those probability distributions as we want a predictor about whether an agent \textit{is going to} switch. With this idea on hand, we can use the aggregation of ex-ante net gains over the population as a disequilibrium index, i.e., a Lyapunov function, to tell whether a population-level dynamic converges to any equilibrium or not: the dynamic should stop at a stationary state if and only if the aggregate gain is zero. This is our first main result (\Cref{clm:GZero,clm:GLyapunov} in \Cref{sec:ExAnteNetGains}).

Dynamic stability of an equilibrium is implied from a monotone decrease in the aggregate net gain toward zero. So, what guarantees monotone decrease? Economists have developed notions of static stability as a property of equilibria or games, without specifying dynamics explicitly. For instance, evolutionary stable states (ESS) is an example of a local static stability condition for an equilibrium and a stable game is regarded as a game with global static stability (H\&S). Mathematically, static stability is summarized by negative semidefiniteness of the Jacobian of a payoff function, which defines static stability in this paper (\Cref{dfn:StatStbl}). Economists conjecture that static stability should imply dynamic stability in economic dynamics: as long as agents respond somehow positively to incentives, then the negative feedback in payoffs guides them back to an equilibrium. In such games, incentives for changing strategies should diminish as agents move toward equilibria.

Despite such a commonly held belief in the link between static and dynamic stability, evolutionary game theory has obtained both positive results (summarized in H\&S) and negative results (e.g. \citealt{Friedman91Ecta_EvolGamesEcon}) about the link. Our notion of net gains can explain \textit{both}. On the one hand, convergence is guaranteed for the best response dynamic (and its variants), pairwise comparison dynamics, and excess payoff dynamics (H\&S). About these positive cases, we propose a general property that encompasses these dynamics. That is \textit{full rationalizability} of myopic decisions.\footnote{One might simply call it rationalizability. The modifier ``cost-benefit'' is added to distinguish it from rationalizability of a strategy as \cite{Pearce_84Ecta_Rationalizable,Bernheim_84Ecta_Rationalizable}. Here we mean by rationalizability that an agent's choice can be explained as a result of some \textit{optimization}, possibly distorted by switching costs or constraints.} When the available strategy set depends on the current strategy that an agent has been taking, a rational agent may not simply switch to the strategy that yields the highest payoff among all the currently available strategies, but to a suboptimal strategy that, however, gives a greater accessibility to a much better strategy than those currently available strategies in the next revision opportunity. This contradicts with the assumption that agents' choices of new strategies are based solely on myopic current payoffs, commonly assumed in evolutionary dynamics. Full rationaliability of myopic decisions will be presented as a condition on \textit{availability }of strategies to avoid such a contradiction (\Cref{clm:g_pi} in \Cref{sec:RatMyop}). Then, in our second main result (\Cref{clm:DG,clm:StatDynStbl} in \Cref{sec:StatStbl_NetGains}), we verify that full rationalizability, coupled with static stability of a game, indeed implies monotone decrease in the aggregate ex-ante net gain and thus dynamic stability of equlibria. 

On the other hand, imitative dynamics such as replicator dynamics do not guarantee dynamic stability. Imitation through sampling of other agents' strategies creates positive feedback to deviation from an equilibrium through changes in the sampling probability of each strategy. This prevents those dynamics from convergence despite negative feedback through changes in payoffs by static stability. The concept of net gains gives a quantitative explanation; we show that the positive feedback through sampling indeed refuels the net gains (\Cref{clm:ReplRefuelG} in \Cref{sec:ImitDyn}). In particular, zero-sum games show the most distinguishable contrast between positive cases and negative cases. While experiments as in \citet{ONeil_87_PNAS_TestZeroSum} show robust empirical supports for convergence to an equilibrium when real human beings play zero-sum games, it is theoretically known that imitative dynamics do not converge. Our results explain this gap. 

The net gain approach has advantages in applications, not only for applicability to a wide variety of dynamics and games whether or not a dynamic converges, but also for \textit{not} requiring us (analysts of games) to have any prior knowledge about a limit state or convergence. This should be useful for prediction of convergence when a game has multiple equilibria. It is indeed in such a situation when we may want to rely on simulations of evolutionary dynamics to find equilibria numerically. This is a technical advantage of our concept of net gains over the known Lyapunov functions in the preceding literature on equilibrium stability in evolutionary and learning dynamics, such as the potential of a game \citep{Sandholm01JET_Potential} and geometric distance/divergence from a limit state \citep{MerikopoulousSandholm18JET_RiemannianGameDyn}. Furthermore, the idea of net gains comes from economic intuition, not from mathematical specification. While here we confine to strategic forms of games, it is easy to extend the approach here to more complicated situations; for example, \cite{SawaZusai_multitaskBRD} applies this approach to multitasking environments, i.e., concurrent plays of multiple games.  

\Cref{sec:Prelim} summarizes basic concepts and notation about evolutionary dynamics in population games. The basic framework of cost-benefit rationalizable dynamics is presented in \Cref{sec:Basic}. We define the aggregate net gain and utilize it to obtain the two main results in \Cref{sec:MainResults}. In \Cref{sec:BehindGain}, we examine economic meanings of full rationalizability and causes of instability in imitative dynamics. \Cref{sec:Ext} collects extensions beyond the basic set-up. In \Cref{sec:dis}, we compare our approach with others in the preceding literature. \Cref{sec:Concl} provides a quick summary of our approach and its implications. The proofs and details of side issues are given in Appendix.

\section{Preliminaries}\label{sec:Prelim}

\subsubsection*{Notation} A vector in a bold font like $\v$ is a column vector, while one with an arrow over the letter like $\vec{v}$ is a row vector. We omit the transpose when we write a column vector on the text. $\1$ is a column vector $(1,1,\ldots,1)$. Note that $\1\cdot\z=\sum_{i=1}^n z_i$ for an arbitrary column vector $\z=(z_i)_{i=1}^n\in\R^n$. For a finite set $\cA=\{1,\ldots,A\}$, we define $\Delta^\cA$ as $\Delta^\cA \coloneqq \{\x=(x_a)_{a\in\cA}\in[0,1]^{\sharp\cA} \mid \1\cdot\x=1 \}$, i.e., the set of all probability distributions on $\cA$. Let $\R_+ \coloneqq [0,\infty)$ and $\R_{++} \coloneqq (0,\infty)$; $\R^\cA \coloneqq \R^{\sharp\cA}$, $\R^\cA_+ \coloneqq \R^{\sharp\cA}_+$, and $\R^\cA_0 \coloneqq\{ \z\in\R^\cA \mid \1\cdot\z=0 \}$.

Consider an $A$-dimensional real space, each of whose coordinate is labeled with either one element of $\cA=\{1,\ldots,A\}$. For a set $S\subset \cA$, we denote by $\Delta^\cA(S)$ the set of probability vectors whose support is contained in $S$: i.e., $\Delta^\cA(S)\coloneqq \left\{ \x\in\R^\cA_+ \mid \sum_{k\in S}x_k=1 \text{ and }x_l=0 \text{ for any }l\in\cA\setminus S \right\}$. For each $a\in\cA$, define a unit vector $\e_a=(e_{ab})_{b\in\cA} \in\Delta^\cA$ such that $e_{aa}=1$ and $e_{ab}=0$ for all $b\ne a$. This represents a population state where every agent takes action $a$. 

In a vector space $\cZ$, with a set $S\subset\cZ$ , an element $\c\in\cZ$ and a scalar $k\in\R$, we define set $ k(S+\c)\coloneqq \{ k(\z+\c) \in\cZ \mid \z \in S \}$; with a scalar set $K\subset\R$, set $KS$ is defined as $KS\coloneqq \{k\z\in\cZ \mid k\in K, \z\in S\}$. For sets $S_1,S_2\subset\cZ$, we define set $S_1+S_2\coloneqq \{\z_1+\z_2\in\cZ \mid \z_1\in S_1 \text{ and } \z_2\in S_2 \}$.

$[\cdot]_+$ indicates the non-negative part of the argument: i.e, $[z]_+=\max\{z,0\}$. $\E_Q$ is the expected value operator: $\E_Q f(q)=\int_\R f(q)d\P_Q(q)$ for integrable function $f:\R\to\R$. 

\subsection{Population game $\F$}
We consider a finite-action population game played in the society of continuously many homogeneous agents, as in \citet[Ch.2]{SandholmPopText}. The society consists of a unit mass of agents. In this basic set-up, they are homogeneous in the sense that they have the same action set and the same payoff function; see \Cref{sec:hetero} for an extension to a heterogeneous setting. 

Specifically, each agent chooses an action $a$ from a finite set $\cA$. Denote by $x_a\in[0,1]$ the mass of action-$a$ players in the society. The \textbf{social state }is represented by column vector $\x\coloneqq(x_a)_{a\in\cA}$ in $\Delta^\cA$. The tangent space of $\Delta^\cA$ is $\R^\cA_0.$ 
Let $F_a(\x)\in\R$ be the payoff from action $a\in \cA$ in social state $\x\in\Delta^\cA$. Define payoff function $\F:\Delta^\cA\to\R^\cA$ by $ \F(\x)\coloneqq (F_a(\x))_{a\in\cA}$ for each $\x\in\Delta^\cA.$ Throughout this paper, we assume that $\F$ is continuously differentiable. Let $F_\ast(\x)\coloneqq \max_{a\in\cA} F_a(\x)$. A population game is characterized by $(\cA,\F)$, or shortly by $\F$. Social state $\x^\ast\in \Delta^\cA$ is a {\bf Nash equilibrium} if $F_a(\x^\ast)=F_\ast(\x^\ast)$ whenever $x_a^\ast>0$. %
Denote by $\NE(\F)$ the set of Nash equilibria in  population game $\F$.

\subsection{Evolutionary dynamic $\V$}\label{sec:EvolDynV}
For now we simply define \textbf{an evolutionary dynamic } $\V:\Delta^\cA\times\R^\cA \to \R^\cA_+$ as a (set-valued) differential equation such as $\dot\x\in \V(\x,\bpi)$; so, it is a deterministic dynamic over a continuous time horizon. $\V$ is set-valued, i.e., a differential inclusion. This is to allow for multiple transition vectors, which is typically a case in the best response dynamic (BRD) due to multiplicity of best responses. Behind this population-level dynamic, each agent occasionally receives an opportunity to revise the action and then decide on whether to switch the action and, if so, what action to switch. This decision is supposed to follow a certain rule, called a \textbf{revision protocol}; a variety of evolutionary dynamics comes from a variety of revision protocols \citep{SandholmPopText}. We will put them under one unifying framework in the next section. 

Once we set up both game $\F:\x\mapsto\bpi$ and evolutionary dynamic $\V:(\x,\bpi)\mapsto\dot\x$, we obtain an autonomous dynamic system $\V^\F:\x\mapsto\dot\x$ to completely determine the transition $\dot\x$ of social state from $\x$ alone: 
$$\dot\x \in \V(\x,\F(\x))=:\V^\F(\x).$$
We call $\V^\F$ the \textbf{combined dynamic}. We say that $\x$ is a \textbf{rest point }in $\V^\F$ if $\0\in \V^\F(\x)$; let $\RP(\V^\F)$ be the set of rest points. 

Dynamic stability of an equilibrium in an autonomous dynamic system simply means asymptotic stability, i.e., convergence (attracting) and no farther escape (Lyapunov stability); a formal explanation is given in \Cref{apdx:DynMath}. In mathematical study of dynamical systems, a Lyapunov function is a common tool to prove or disprove dynamic stability. Say, we have a scalar-valued continuous function $W^\F:\Delta^\cA\to \R$ whose zeros coincide with rest points of $\V^\F$, i.e., $W^\F(\x)=0 \ \Leftrightarrow \x\in\RP(\V^\F)$. If $W^\F(\x^t)$ converges to zero on every solution trajectory $\{\x^t\}_{t\in\R_+}$ of $\dot\x^t\in\V^\F(\x^t)$ from any initial point $\x^0$, then the set of rest points is asymptotically stable. If $W^\F(\x^t)$ diverges from zero from any initial point around a rest point $\x^\ast$, then this rest point is unstable. Such a function is a \textbf{Lyapunov function} (see \Cref{thm:Lyapunov_DI} in \Cref{apdx:DynMath}). The study of multi-dimensional dynamic systems boils down to tracking the value of this function $W^\F$. 

\section{Basic framework}\label{sec:Basic}

\subsection{Rationalizing bounded rationality}
To measure net gains, we rationalize bounded rationality in an evolutionary dynamic. Inertia or reluctance to switch is rationalized by introducing \textit{switching costs}. Switches to suboptimal strategies are rationalized by constraints on available strategies so we say that the optimal strategy was \textit{not available}. To explain continuity, we allow for randomness in switching costs and available strategy sets. Then, we reconstruct a revision protocol from optimization under a stochastic switching cost subject to a stochastic available strategy set.\footnote{Stochastic switching cost alone is considered in \cite{Kuzmics10DomStr_StochBRD} on stochastic evolution and in \citet{ZusaiTBRD} on deterministic evolution (tempered BRD); stochastic restriction to available strategies set alone is considered in \cite{BravoFaure15SIAMJCtrl_ReinforceLearn_RestrAction}.} 

We make it formal. Over a continuous time horizon, an agent occasionally revises the action in the following steps.
\begin{description}
	\item[Step 0] According to a Poisson process with arrival rate 1, a revision opportunity arrives at an agent. Say, the agent has been taking action $a$ until this moment of time.
	\item[Step 1] Upon receipt of a revision opportunity, the agent draws \textit{a set of available actions} $\prA\subset\cA$ from probability distribution $\P_{\cA }$ over the power set of $\cA$. Observing payoff vector $\bpi$, the agent finds the best \textit{available} action among actions in $\prA$; let $b_\ast(\bpi,\prA) \coloneqq  \argmax_{b\in \prA}  \pi_b$ be the set of these best available actions in $\prA$, and $\pi_\ast(\prA)\coloneqq  \max_{b\in \prA}  \pi_b$ be the payoff from these actions.
	\item[Step 2]  The agent also draws the value of \textit{a switching cost} $q\in\R_+$ from probability distribution $\P_Q$ over $\R_+$; it must be paid if an agent switches the action. Thus, the agent chooses to switch to the best available action (if more than one, pick any) by paying $q$, if  $\pi_\ast(\prA)-q$ is greater than the current payoff $\pi_a$; the agent chooses not and keeps taking action $a$, if $\pi_\ast(\prA)-q<\pi_a$.   
\end{description}
These steps are summarized into the maximization problem below. 
\begin{equation}
	\max \left\{ \pi_a, \max_{b\in \prA}  \pi_b-q \right\} \label{eq:OptProtocol}
\end{equation}
To trigger a switch, the \textbf{gross gain} from the switch $\breve\pi_\ast(\prA):=\pi_\ast(\prA)-\pi_a$ should exceed the switching cost $q$. We call $\breve\pi_\ast(\prA)-q$ the \textbf{net gain}.

A variety of revision protocols reduces to a variety of $\P_\cA$ and $\P_Q$. A revising agent may stick with the current action despite a positive gross gain; this is now rationalized as the gross gain was smaller than $q$. The new action $b$ may not be the exact optimal action $b_\ast(\bpi)\coloneqq b_\ast(\bpi,\cA)$ among all the actions; it is now rationalized as any of $b_\ast(\bpi)$ was out of the available action set $\prA$ and the new action $b$ was still the best among the available actions in $\prA$, i.e., $b\in b_\ast(\bpi,\prA)$. Note that yet we require new action $b$ not to yield a smaller payoff than the old action, i.e., $b$ to be a better reply to $\bpi$ than $a$.

We do not impose any tie-breaking rule when there are multiple best available actions (i.e., $b_\ast(\bpi,\prA)$ is not a singleton) or when an agent is indifferent between switching and no switching (i.e., $\pi_\ast(\prA)-q=\pi_a$). In these cases, an agent's choice may be probabilistic and the probability may exhibit some indeterminacy. From an ex-ante viewpoint before $q$ is drawn but after $\prA$ is drawn, the probability of choosing a switch can take any value in the range $\cQ(\breve\pi_{a\ast}(\prA))$, where $\cQ(\bar q)\coloneqq [Q_-(\bar q), Q(\bar q)]$ with $Q_-(\bar q)\coloneqq\P_Q(\{q<\bar q\}))$ and $Q(\bar q)\coloneqq\P_Q(\{q\le \bar q\})$.\footnote{If $Q$ is continuous at $\bar q$, i.e., $Q_-(\bar q)=Q(\bar q)$, then $\cQ(\bar q)$ reduces to a singleton $\{Q(\bar q)\}$. While we distinguish $\cQ$ from $Q$ for mathematical completeness in case $\P_Q$ is not continuous, readers do not have to be bothered when reading the main text; when it matters anyhow, it is cautioned explicitly.} Conditional on choosing to switch, the probability distribution of a new action takes any in $\Delta^\cA(b_\ast(\bpi,A'))$.

For simplicity in exposition, we assume that availability of actions $\P_\cA$ does not depend on the current action until \Cref{sec:RatMyop}. It will be relaxed and some theorems do not need this assumption, so we make it explicit below.\footnote{We allow $\prA$ to include the current action, though the current action is always available as a status quo without paying a switching cost in this basic framework.} We allow $\P_\cA$ to depend on $\x$ to cover imitative dynamics. Denote by   $\P_\cA(\prA;\x)$ the probability with which $\prA$ is the available action set given $\x$. We assume that $\P_Q$ has no singular continuous part\footnote{That is, $\P_Q$ is a convex combination of an absolute continuous measure and a discrete measure.} and that $\prA$ and $q$ are independent of each other. 
\begin{assmpA}\label{assmp:AIndp_a}
	For any $\prA\subset\cA$, the probability $\P_\cA(\prA;\x)$ does not depend on a revising agent's current action $a\in\cA$.
\end{assmpA}

\subsection{Aggregation to an evolutionary dynamic}

We obtain a population-level evolutionary dynamic by aggregating individual revision processes over the entire population. We say that an evolutionary dynamic $\V$ is \textbf{cost-benefit rationalizable} if $\V$ is constructed from a constrained optimization protocol \eqref{eq:OptProtocol} with some $\P_\cA$ and $\P_Q$.  The basic formulation in this section captures not only the standard best response dynamic (BRD) but also tempered BRDs and pairwise comparison dynamics; see \Cref{fig:majorDyn}. It also allows for new evolutionary dynamics, such as \textit{ordinal better-reply dynamics}, where switching rates depend on the \textit{ordinal }payoff ranking.\footnote{See \Cref{fig:majorDyn} for $\P_\cA$ to rationalize this dynamic.  In $\P_\cA$ for this dynamic, each action becomes available with an equal probability, independently of availability of other actions. The constant $c$ is for normalization so that the sum of $\P_\cA(\prA)$ over all non-empty sets $\prA$ equals to 1. Bill Sandholm suggested that this $\P_{Aa}$ implies ordinal better-reply dynamics.} Our \textit{approach }covers excess payoff dynamics and smooth BRDs with a little modification of the formulation, particularly in the definition of a ``status-quo'' strategy for an agent, to be discussed in \Cref{sec:Modified}.

Below we present the formula of $\V$ in a cost-benefit rationalizable dynamic. (A reader who is familiar with construction of a population-level evolutionary dynamic from a revision protocol as in \citet[Chapter 4]{SandholmPopText}, may skip the rest of this subsection.) Despite randomness in revision opportunities, the assumption of a continuous population allows us to formulate the transition of the social state $\dot\x$ as a deterministic dynamic \citep{RothSandholm_StochApproxDI}.  
 
Evolutionary dynamic $\V:\Delta^\cA\times\R^\cA \rightrightarrows\R^\cA_0$ is constructed from $\P_\cA$ and $\P_Q$ as 
\begin{equation}\label{eq:TransitSet}
\dot\x\in \V(\x,\bpi):= \sum_{a\in\cA} x_a \underbrace{\sum_{\prA\subset \cA} \P_{\cA a}(\prA;\x) \cQ(\breve\pi_{a\ast}(\prA)) \left\{ \Delta^\cA( b_\ast(\bpi,\prA))-\e_a \right\} }_{\displaystyle \V_a(\bpi)} .
\end{equation}
The right-hand side is interpreted as follows. Consider agents who receive revision opportunities in infinitesimal time period $[t,t+dt)$. The mass of those agents is $1\times dt$ since the arrival rate of revision opportunities is $1$. Among them, $x_a$ is the proportion of those who have been taking action $a$ until the receipt of revision opportunities in step 0. In step 1, for each $\prA\subset\cA$, $\P_{\cA a}(\prA)$ is the proportion of agents who face available strategy set $\prA$ among those. In step 2, they choose to abandon action $a$ with any probability in the range $\cQ(\breve\pi_{a\ast}(\prA))$. As they switch to any of the best available strategies in $b_\ast(\bpi,\prA)$, the distribution of their new actions can be any probability distribution in $\Delta^\cA( b_\ast(\bpi,\prA))$ while that of their old actions is $\e_a$ since they were all taking $a$. The transition of the social state $d\x$ in this time period is obtained by aggregating these changes over all $a\in\cA$ and $\prA\subset\cA$.  

\newcommand{\highlight}[1]{{\color{red}#1}}
\newcommand{\dynEntry}[4]{
	\begin{tcolorbox}[arc=3mm,boxrule=0.5mm,boxsep=0.5mm,
		fonttitle=\sffamily\bfseries,fontupper=\small,
		coltitle=black,	colbacktitle=white,colback=white,
		title=#1]
		#2
		#3
		
		#4
\end{tcolorbox}}
\newcommand{\dynEntryFull}[4]{
	\begin{tcolorbox}[arc=3mm,boxrule=0.5mm,boxsep=0.5mm,
		fonttitle=\sffamily\bfseries,fontupper=\small,
		coltitle=white,	colbacktitle=black,colback=white,
		title=#1]
		#2
		#3
		
		#4
\end{tcolorbox}}
\newcommand{\dynEntryFullMod}[4]{
	\begin{tcolorbox}[enhanced,arc=3mm,boxrule=0.5mm,boxsep=0.5mm,
		borderline={0.3mm}{0mm}{black},borderline={0.3mm}{0mm}{white,dashed},
		fonttitle=\sffamily\bfseries,fontupper=\small,
		coltitle=white,	colbacktitle=black,colback=white,
		title=#1]
		#2
		#3
		
		#4
\end{tcolorbox}}

\begin{figure}[p!]
	
	\dynEntryFull{Tempered BRDs \citep{ZusaiTBRD}}
	{Switch to $b_\ast(\bpi)$ with probability $Q(\pi_\ast-\pi_a)$.}
	{$\Leftarrow$ \quad $\P_\cA(\cA)=1$ with arbitrary $Q$.}
	{%
		\dynEntryFull{Standard BRD \citep{Hofbauer95BRD}}
		{Always switch to $b_\ast(\bpi)$ regardless of $\pi_\ast-\pi_a$.}
		{$\Leftarrow$ \quad $\P_\cA(A)=1$. $Q(\bar q)=1$ for any $\bar q>0$.}
		{}
	}
	\dynEntryFull{Pairwise comparison dynamics (H\& S)}
	{Sample each action $b\in\cA$ with an equal probability, and switch to it with probability $Q(\pi_b -\pi_a)$.}
	{$\Leftarrow$ \quad $\P_\cA(\{b\})=1/A$ for each $b\in\cA$, with arbitrary $Q$.}
	{%
		\dynEntryFull{Smith dynamic \citep{Smith84TranspSci_Stbl_Dyn_TrafficAssignmt}}
		{Pairwise comparison dynamic with $Q(\pi_b -\pi_a)=[\pi_b-\pi_a]_+.$}
		{$\Leftarrow$ \quad $\P_\cA(\{b\})=1/A$ for each $b\in\cA$. $Q(\bar q)=[\bar q]_+$.}
		{}			
	}
	\dynEntryFull{Ordinal better-reply dynamics (new in this paper)}
	{Switch to the $i$-th best action with probability $c\bar{p}(1-\bar{p})^{i-1}$ as long as it is better than $a$. ($\bar{p}\in(0,1)$ and $c=1/(1-(1-\bar{p})^{A-1})$.)}
	{$\Leftarrow$ \quad $\P_{\cA}(\prA)=c\bar{p}^{\sharp\prA}$ for each $\prA\ne\emptyset$. $Q(\bar q)=1$ for any $\bar q>0$.}
	{}
	
	\dynEntry{Imitative dynamics }
	{Sample each action $b\in\cA$ with probability $x_b$, 	and switch to it with probability $Q(\pi_b -\pi_a)$.}
	{$\Leftarrow$ \quad $\P_\cA(\{b\})=x_b$ for each $b\in\cA$, with arbitrary $Q$.}
	{%
		\dynEntry{Replicator dynamic \citep{Smith84TranspSci_Stbl_Dyn_TrafficAssignmt}}
		{Imitative dynamic with $Q(\pi_b -\pi_a)=[\pi_b-\pi_a]_+.$}
		{$\Leftarrow$ \quad $\P_\cA(\{b\})=x_b$ for each $b\in\cA$. $Q(\bar q)=[\bar q]_+$.}
		{}			
	}
	
	\dynEntryFullMod{(Separable) excess payoff dynamics (H\&S)}
	{Sample an action $b\in\cA$ with probability $1/A$ and switch to it with probability increasing its relative payoff $\hat\pi_b$.}
	{$\Leftarrow$ \quad $\P_{\cA}(\{b\})=1/A$ for each $b$ with any $Q$}
	{
		\dynEntryFullMod{Brown-von Neuman-Nash dynamic \citep{Hofbauer00Sel_NashBrownToMaynardSmith}}
		{Excess payoff dynamic with $Q(\hat\pi)=[\hat\pi]_+$.}
		{$\Leftarrow$ \quad $\P_{\cA}(\{b\})=1/A$ for each $b$ and $Q(\bar q)=[\bar q]_+$}
		{}
	}
	
	\caption{Major dynamics captured in cost-benefit rationalizability. The dynamics in black title boxes are fully rationalizable. A box with dashed borderline means that the dynamic needs a modification of the framework as in \Cref{apdx:Modified}. A box containing another box suggests that the class of dynamics represented in a larger box includes the dynamics in a smaller box. In each box, a revision protocol is described first and followed by the pair of $\P_{\cA a}$ and $Q$ to rationalize it; $a$ denotes the action that a revising agent has been taking until receiving a revision opportunity and $\bpi$ is the payoff vector at this moment of time. Each citation refers to a paper that derives dynamic stability from static one (stable games or ESS) under the dynamic in the box.
	}\label{fig:majorDyn}
\end{figure}

\section{Main results \label{sec:MainResults}}

\subsection{Ex-ante net gains\label{sec:ExAnteNetGains}}
\subsubsection*{Definition and the link with rest points}
While a gross gain (without considering $\prA$ or $q$) does not explain an agent's decision on whether and to what to switch from current action, say action $a$, we can now say that this decision is made to maximize the net gain, i.e.,
$$ \max \bigg\{ 0, \underbrace{\max_{b\in \prA}  \pi_b-\pi_a}_{\displaystyle \breve\pi_{a\ast}(\prA)} -q \bigg\}. $$ 
As this maximal value is shortly rewritten as $[\breve\pi_{a\ast}(\prA)-q]_+$, we define the \textbf{individual ex-ante (first-order) net gain} for an action-$a$ player as \footnote{Note that $g_{a\ast}[\bpi]$ is well defined: that is, the value of the RHS of \eqref{eq:g_dfn} is uniquely determined without choosing a particular value of a switching probability from $\cQ(\breve\pi_{a\ast}(\prA))$ or assuming continuity of $\P_Q$, regardless of whether or not to switch when the net gain is exactly zero.} 
\begin{equation}
	g_{a\ast}(\x,\bpi)\coloneqq   \sum_{\prA \subset \cA } \P_\cA(\prA; \x) \E_Q[\breve\pi_{a\ast}(\prA)-q]_+. \label{eq:g_dfn}
\end{equation}
This represents the remaining room to improve the payoff for a player who currently takes action $a$. The net gain $\breve\pi_{a\ast}(\prA)-q$ from a switch is counted into the expected value $\E_Q[\breve\pi_{a\ast}(\prA)-q]_+$ if and only if the agent chooses a  switch, i.e., the net gain is positive. 

Note that $\x$ affects $g_{a\ast}(\x,\bpi)$ only through $\P_{\cA a}(\prA; \x)$, i.e., only through changes in availability of strategies. Thus, $g_{a\ast}$ does not depend on $\x$, if availability of actions $\P_\cA$ does not depend on the current social state $\x$, as in the assumption below.
\begin{assmpA}\label{assmp:AIndp_x}
	For any $\prA\subset\cA$, $\P_\cA(\prA;\x)$ does not depend on $\x$.
\end{assmpA}
Imitative dynamics including the replicator dynamic do not satisfy this assumption. It is known that imitative dynamics fail to assure asymptotic stability of rest points in zero-sum games, unlike non-imitative (payoff-based) dynamics. We will relate this instability with violation of \Cref{assmp:AIndp_x} in \Cref{sec:ImitDyn} through a close investigation of causes to change net gains.

To analyze the population dynamic of $\x$, we define the \textbf{aggregate (first-order) net gain function }
$G:\Delta^\cA\times\R^\cA \to\R$ as 
$$ G(\x,\bpi) \coloneqq \x\cdot \g_\ast(\x,\bpi) =\sum_{a\in\cA} x_a g_{a\ast}(\x,\bpi),$$ 
where $\g_\ast(\x,\bpi) \coloneqq ( g_{a\ast}(\x,\bpi) )_{a\in\cA} \in\R^\cA$. As we combine $\F$ and $\V$ to obtain an autonomous dynamic $\V^\F$, we embed $\F$ to $G$ and define $G^\F:\Delta^\cA \to\R$ by $G^\F(\x) \coloneqq G(\x,\F(\x))$.

Now we establish a foundation to use the aggregate net gain $G^\F$ as a Lyapunov function. 
\begin{thm}\label{clm:GZero}
	Consider a cost-benefit rationalizable dynamic $\V$. Then, the aggregate net gain function $G:\Delta^\cA\times\R^\cA \to\R$ satisfies the following property. 
	\begin{description}
		\item[G0-a)\label{prop:GZero}] i) $G(\x,\bpi)\ge 0$ for any $\x\in\Delta^\cA,\bpi\in\R^\cA;$ ii) $G(\x,\bpi)=0$ if and only if $ \0\in\V(\x,\bpi) $.
	\end{description}
\end{thm}
Property \ref{prop:GZero} suggests that convergence of $\x$ to rest points in $\V^\F$ is equivalent to convergence of $G^\F(\x)$ to zero.  This coincidence is generally guaranteed, even if $\RP(\V^\F)\neq \NE(\F)$, because of counting switching cost $q$ in $G^\F$. This is the first stark difference from the gross gain, to be discussed in \Cref{sec:GrossGain}.

\begin{cor}\label{clm:GLyapunov}
	Consider a cost-benefit rationalizable dynamic $\V$ in a population game $\F$. Then, $G^\F$ serves as a Lyapunov function: convergence in $G^\F$ to zero implies convergence of $\x$ to rest points and divergence in $G^\F$ from zero implies instability of a rest point.
\end{cor}

\subsubsection*{Three factors to change $G^\F$}
According to \Cref{clm:GLyapunov}, monotone \textit{decrease }in $G^\F$ to zero implies asymptotic stability of rest points  under $\V^\F$, while monotone \textit{increase }in $G^\F$ from zero implies instability. Thus, we can find causes of stability or instability by looking into effects of change $\dot\x$ in the social state on aggregate net gain $G^\F$. The instantaneous change $\dot G^\F$ in the aggregate net gain is decomposed into three factors:\footnote{For differentiable function $\f=(f_i)_{i=1}^m:\R^n\to\R^m$ with input variable $\z=(z_j)_{j=1}^n$, we denote the Jacobian of $\f$ by $D\f$ or $d\f / d\z$; $\tpdif{f_i}{z_j}$ locates on the $i$-th row and the $j$-th column of  $d\f / d\z$.}
$$ \dot G^\F(\x) = \underbrace{\x \cdot \pdif{\g_\ast}{\bpi}(\x,\F(\x)) \overbrace{\pdif{\F}{\x}(\x) \dot\x}^{\dot\bpi} }_{\substack{\displaystyle\text{Payoff effect}\\ {\small\displaystyle \pdif{G}{\bpi} \dot\bpi}} }  + \underbrace{\underbrace{\x \cdot \pdif{\g_\ast}{\x}(\x,\F(\x)) \dot\x}_{\displaystyle\text{Availability effect}}  + \underbrace{\dot\x \cdot \g_\ast(\x,\F(\x))}_{\displaystyle\text{Switching effect}} }_{\small\displaystyle \pdif{G}{\x} \dot\x}.$$
The \textbf{payoff effect }aggregates the changes in individual agents' net gains $\g_\ast$ \textit{caused by the change in payoff vector }$\dot\bpi = D\F(\x) \dot\x$ over all agents; the payoff change $\dot\bpi$ affects a revising agent's assessment of ex-post net gain $[\breve\pi_\ast(\prA)-q ]_+$  and thus an agent's individual ex-ante net gain $\g_\ast$, which are aggregated to $G^\F$. The \textbf{availability effect }aggregates the indirect effects of changes in $\x$ on individual net gains $\g_\ast$ \textit{through changing availability of each strategy}; note that, provided that $\bpi$ was fixed, $\x$ should affect $\g_\ast(\x,\bpi)$ only through $\P_\cA(\prA;\x)$. To understand the \textbf{switching effect}, notice that $x_a$ represents not only the share of action-$a$ players but also the share of $g_{a\ast}$  in the distribution of individual ex-ante net gains over all the agents and thus $\dot\x$ represents the change in this distribution; a switch from action $a$ to action $b$ changes the switching agent's ex-ante gain from $g_{a \ast}(\x,\bpi)$ to $g_{b \ast}(\x,\bpi)$. Such changes in revising agents' individual ex-ante net gains \textit{by their own switches }are collected in the switching effect.

\Cref{clm:DG} relates these three effects with natures of game $\F$ and dynamic $\V$. 
\begin{thm}\label{clm:DG}
	Consider a cost-benefit rationalizable dynamic $\V$. Then, the aggregate net gain function $G:\Delta^\cA\times\R^\cA \to\R$ satisfies the following property. 
	\begin{description}
		\item[G0-b)\label{prop:GDpi}] (Payoff effect) $G$ is Lipschitz continuous and thus differentiable almost everywhere.  Whenever  $G$ is differentiable at $(\x,\bpi)\in\Delta^\cA\times\R^\cA,$ we have
		$$ \pdif{G}{\bpi}(\x,\bpi) \Delta\bpi = \Delta\x\cdot\Delta\bpi \qquad\text{		for any $\Delta\x\in\V (\x,\bpi), \Delta\bpi\in\R^\cA.$}$$
	\end{description}
	Furthermore, we have the following.
	\begin{description}
		\item[G1 \label{prop:gzNeg}] (Switching effect) Under \Cref{assmp:AIndp_a}, there exists a function $H:\Delta^\cA\times\R^\cA \to\R$ such that 
		\begin{enumerate}[a)]
			\item $H$ is lower semi-continuous; besides, i) $H(\x,\bpi)\le 0$ at any $\x,\bpi$ and ii) $H(\x,\bpi)=0$ if and only if $\0\in \V(\x,\bpi)$;
			\item $ H(\x,\bpi) \equiv \g_\ast(\x,\bpi)\cdot \Delta\x$ for any $\Delta\x\in\V(\x,\bpi)$.
		\end{enumerate}
		\item[G2  \label{prop:gDx}] (Availability effect) Under \Cref{assmp:AIndp_x}, $\g$ does not change with $\x$ alone, i.e., $\g_\ast(\x,\bpi)=\g_\ast(\x',\bpi)$  for any pair of $\x,\x'$ and for any $\bpi$. Hence, we have 
		$$ \pdif{\g_\ast}{\x}(\x,\bpi) = \O. $$
	\end{description}
\end{thm}

To understand property \ref{prop:GDpi}, imagine a change in payoff vector by $\Delta\bpi$ and its effect on the individual ex-ante net gain $g_{a \ast}$ for an action-$a$ player. Given a draw of $\prA$ and $q$, the ex-post net gain $\breve\pi_{a\ast}(\prA)-q$ changes by $\Delta\pi_\ast(\prA)-\Delta\pi_a$ regardless of $q$.\footnote{$\Delta\pi_\ast(\prA)$ is the change in the payoff from the best action $b_\ast(\bpi,\prA)$, i.e., $\Delta\pi_\ast(\prA)=\Delta\pi_b$ of an action $b\in b_\ast(\bpi,\prA)$. The derivative $\tpdif{g_{a\ast}}{\bpi}$ is well defined (when $g_a\ast$ is differentiable) regardless of multiplicity of the best actions, thanks to Dankin's envelop theorem; see the proof of property \ref{prop:gDpi} in \Cref{apdx:Proofs}.}  The resulting change in $g_{a\ast}(\bpi',\x)$ is linearly approximated as\footnote{To obtain this equation mathematically, notice that the derivative of $\E_Q[\breve\pi-q]_+=\int_{q=0}^{\breve\pi} (\breve\pi-q) d\P_Q(q)$ with respect to $\breve\pi$ is $\int_{q=0}^{\breve\pi} 1\cdot d\P_Q(q)=\P_Q(q\le \breve\pi)=Q(\breve\pi)$.}
\begin{align*}
	\pdif{g_{a\ast}}{\bpi} (\x,\bpi) \Delta\bpi
	& = \sum_{\prA \subset \cA } \P_\cA (\prA;\x) Q(\breve\pi_{a\ast}(\prA)) \cdot (\Delta\pi_\ast(\prA)-\Delta\pi_a) \\
	& = \z_a\cdot\Delta\bpi \qquad\text{ with any $\z_a\in\V_a (\x,\bpi) $}.
\end{align*}
This is the expected sum of a change in the net gain $\Delta\pi_\ast(\prA)-\Delta\pi_a$, weighted with the probability $\P_\cA(\prA) Q(\breve\pi_{a\ast}(\prA))$ that a switch from $a$ is chosen, given available action set $\prA$; this weighted sum is simply $\z_a\cdot\Delta\bpi$.\footnote{When $g_{a\ast}$ is differentiable at $\bpi$, then $\P_Q$ must be continuous at $\breve\pi_{a\ast}(\prA)$ of all $\prA\subset\cA$ for all $a\in\cA$.} Aggregation of this term over all agents yields property \ref{prop:GDpi}.  It is crucial for property \ref{prop:gDpi} that \textit{zero }net gain sets a threshold between switching to a new action and keeping the current action. This marks the second stark contrast with a gross gain, as we will argue in \Cref{sec:GrossGain}. 

Property \ref{prop:gzNeg} means a negative correlation between the ex-ante net gain of each action and the change in the mass of players of the action. We will examine the logic behind this property in  \Cref{sec:RatMyop} where we relax \Cref{assmp:AIndp_a}. Property \ref{prop:gDx} is immediate from \eqref{eq:g_dfn}. Notice that, according to properties \ref{prop:GDpi} and \ref{prop:gzNeg}, $\Delta\x\cdot \Delta\bpi$ and $\g_\ast(\x,\bpi)$ are independent of the choice of a transition vector $\Delta\x$ from $\V (\x,\bpi)$, since the left hand sides $\tpdif{G}{\bpi} (\x,\bpi) \Delta\bpi$ and $H(\x,\bpi)$ in these properties are so.
b
\subsection{Static stability and net gains\label{sec:StatStbl_NetGains}}
\subsubsection*{Static stability}
Now we search for conditions to guarantee monotone decrease in $G^\F$ for dynamic stability of dynamic $\V^\F$. Under \Cref{assmp:AIndp_a} implies a negative switching effect on  $\dot G^\F$ and \Cref{assmp:AIndp_x} implies a zero availability effect. Remaining is the payoff effect, which is determined by game $\F$. It is non-positive if $D\F$ is negative semidefinite; then we have negative $\dot G^\F$. 

Negative semidefiniteness of $D\F$ is indeed a mathematical definition of the static notion of stability in economics. Without explicitly formulating how $\x$ adjusts to $\bpi$, conventional economists believe that the economy should return to an equilibrium as long as there is negative correlation between deviations in action distribution $\Delta\x\in \R^\cA_0$ and the resulting changes in payoff vector $\Delta\bpi.$ With a linear approximation of the latter by $\Delta\bpi\approx D\F(\x)\Delta\x$, the negative correlation is expressed as
\begin{equation}
	\Delta\x \cdot D\F(\x)\Delta\x  \le 0 \qquad \text{ for any $\Delta\x\in \R^\cA_0$}.\label{eq:DF_StatStblCondn}
\end{equation}
This is, of course, nothing but negative semidefiniteness of $D\F$. This condition states that, on average over all actions, actions that have greater shares of players in the population than in the equilibrium should not yield greater payoffs than in the equilibrium, and actions with smaller shares should yield greater profits; H\&S call it \textit{self-defeating externality}. Economists \textit{expect }incentive-driven economic agents to switch from the former group of actions to the latter and thus the society to return to the equilibrium.\footnote{\cite{Milchtaich_19_StaticStbl} attempts to exactly \textit{define }a static stability concept that is universally applicable to different settings (finite player games, population games, etc.).}  Note that \eqref{eq:DF_StatStblCondn} imposes a condition only on $\F$, but not on $\V^\F$. So, we can regard it as \textit{static }stability. 

\begin{dfn}[Static stability]\label{dfn:StatStbl}
	Consider a population game $\F:\Delta^\cA\to\R^\cA$. If $D\F(\x)$ satisfies the negative semideniteness condition \eqref{eq:DF_StatStblCondn} at any $\x\in\Delta^\cA$, we say that the game satisfies \textbf{static stability }globally. If \eqref{eq:DF_StatStblCondn} is satisfied at any $\x$ in some neighborhood of $\x^\ast\in\Delta^\cA$, then we say that the game satisfies static stability locally around $\x^\ast$.
\end{dfn}
We do \textit{not} assert the above static stability is a newly invented notion.\footnote{Recall that a classical market model is ultimately represented by an excess demand function $\d=\D(\p)$ that maps price vector $\p$ to excess demands $\d$; Hick's static stability condition is essentially the negative semidefiniteness of Jacobian $D\D$ of the excess demand function. Thus, the negative semidefiniteness of a Jacobian has been a key element in the notion of static stability in economic theory since John Hicks (\citeyear{Hicks_ValueCap}) proposed it.} Rather, it  summarizes the existing concepts of static stability in game theory. Global static stability of $\F$ is indeed equivalent to $\F$ being a \textit{stable game} (H\&S, Theorem 2.1; also known as contractive games or negative definite games). Static stability is locally satisfied at an interior ESS satisfies local static stability with a strict inequality (unless $\Delta\x\ne\0$) in \eqref{eq:DF_StatStblCondn}. Static stability characterizes a weaker version of ESS, a neutrally stable state (NSS). 

\begin{rmk}
	Note that our static stability is equivalent to monotonicity in comparative statics.\footnote{Our finding (esp., \Cref{clm:StatDynStbl}) in this paper can be regarded as a general proof of the classic correspondence principle in \cite{Samuelson_41_ECTA_StblEqm_CompStat_Dyn} for economic dynamics based on (possibly, boundedly) rational choices.} To see it, consider an exogenous change in payoffs $\btheta\in\R^\cA$ in addition to the original payoff function $\F^0$ so the payoff function $\F:\Delta^\cA\times\R^\cA\to\R^\cA$ is now given by  $\F(\x;\btheta):=\F^0(\x)+\btheta$. Assume the existence of an interior equilibrium $\x^\ast_0$ in the original game $\F^0$; it satisfies $\hat\F(\x^\ast_0;\0)=\0$ where the relative payoff vector  $\hat\F(\x;\btheta)\in\R^\cA$ is defined as $\hat F_a(\x;\btheta)=F_a(\x;\btheta)-\F(\x;\btheta)\cdot\x$ for each $a\in\cA$. This implies that, around $(\btheta,\x)=(\0,\x^\ast_0)$, there is $\x^\ast(\btheta)$ for each $\btheta$ such that $\hat\F(\x^\ast (\btheta);\btheta)=\0$ with differentiable mapping $\x^\ast$ from $\btheta$ to $\x^\ast(\btheta)$. Differentiation of this equation yields
	$$ \Delta\x^\ast\cdot D\F^0(\x^\ast_0) \Delta\x^\ast+\Delta\x^\ast \cdot \Delta\btheta =0. $$
	Therefore, static stability \eqref{eq:DF_StatStblCondn} is equivalent to positive correlation (monotonicity) between exogenous changes in payoff vector $\Delta\btheta$ and changes in equilibrium state $\Delta\x^\ast$, i.e., $\Delta\x^\ast \cdot \Delta\btheta\ge 0$; strictness of the former is equivalent to that of the latter. In particular, strict monotonicity implies that an exogenous increase in the payoff of action $a$, ceteris paribus, causes an increase in the equilibrium share of the action.
\end{rmk}

\subsubsection*{Linking static and dynamic stability through the net gain}
We say an evolutionary dynamic is \textbf{\textit{fully }cost-benefit rationalizable }if it is cost-benefit rationalizable and  $\P_\cA$ satisfies \Cref{assmp:AIndp_a_gen,assmp:AIndp_x_gen}; see \Cref{fig:majorDyn} (black title boxes) for \textit{fully} cost-benefit rationalizable dynamics. According to \Cref{clm:DG}, it implies 
all of properties \ref{prop:GZero}, \ref{prop:GDpi}, \ref{prop:gzNeg} and \ref{prop:gDx} hold. These properties jointly imply\footnote{This equation shows that a fully cost-benefit rationalizable dynamic $\V$ is \textit{$\delta$-passive}, which \citet{FoxShamma13Games} defines as a property of $\V$ alone (before embedding $\F$ to get $\V^\F$) that yields a Lyapunov function of $\V^\F$ when embedded into a stable game $\F$.}  
$$ \pdif{G}{\x}(\bpi,\x) \Delta\x + \pdif{G}{\bpi}(\bpi,\x) \Delta\bpi \le \Delta\x \cdot \Delta\bpi \text{ with a strict inequality if }\0\notin \V(\bpi,\x).$$
Static stability of $\F$ further implies
$$ \dot G^\F(\x) \le \dot\x \cdot D\F(\x) \dot\x \le 0 \text{ with strict inequalities if }\x\notin\RP(\V^\F).$$
Thus, under these assumptions, $G^\F$ decreases toward zero until $\x$ reaches $\RP(\V^\F)$. It suggests that $\RP(\V^\F)$ is asymptotically stable. 

\begin{thm}[Dynamic stability from static stability through net gains]\label{clm:StatDynStbl}
	Consider a fully cost-benefit rationalizable dynamic $\V^\F$. Then, static stability of equilibrium in $\F$ (i.e., negative semidefiniteness of $D\F$) implies  monotone decrease in the aggregate net gain $G^\F$ to zero and thus asymptotic stability under $\V^\F$. In particular, the following holds.
	
	\begin{enumerate}[i)]
		\item Global static stability of $\F$ implies global asymptotic stability of $\RP(\V^\F)$ under the combined dynamic $\V^\F$.
		\item In an arbitrary population game $\F$, let $\x^\ast\in\Delta^\cA$ be an isolated rest point under $\V^\F$. (a) Local static stability of $\F$ at $\x^\ast$ implies (local) asymptotic stability of $\x^\ast$ under the combined dynamic $\V^\F$. (b) Consider a case that $D\F(\x^\ast)$ is not negative semidefinite. Assume that there is indeed a transition vector $\z\in\V^\F(\x^\ast)$ such that $\z\cdot D\F(\x^\ast)\z>0$.\footnotemark ~ Then, $\x^\ast$ is not Lyapunov stable.
	\end{enumerate}
\end{thm}
\footnotetext{Notice this additional assumption. So, this theorem still leaves a room for local stability of $\x^\ast$ if $D\F(\x^\ast)$ is indefinite.}
In short, this theorem tells that \textit{full rationality} of an evolutionary dynamic is sufficient to complete the link between static and dynamic stability. Currently, full rationality may appear to be just a collection of assumptions, i.e., cost-benefit rationalizability and \Cref{assmp:AIndp_a,assmp:AIndp_x}. In \Cref{sec:RatMyop}, we see that these assumptions indeed boil down to a single principle, i.e., rationalizability of myopic payoff-based decision, which underlies \textit{evolutionary }dynamics.

\begin{rmk}
	One might further want to obtain stability of \textit{Nash equilibria}.  Our construction of an optimization-based revision protocol implies that a Nash equilibrium is a stationary state, i.e., $\NE(\F)\subset \RP(\V^\F)$. Furthermore, the following assumptions on $\P_Q$ and $\P_\cA$ imply $\NE(\F)=\RP(\V^\F)$. However, they are not needed for any of the theorems so far and thus \textit{not essential }for dynamic stability.
	
	\begin{assmpQ}\label{assmp:qSupp}  
		$Q(\bar q)>0$ at any $\bar q>0$.%
	\end{assmpQ}
	\begin{assmpA}\label{assmp:AExAnteavailable}
		For any $b\in\cA$, $\P_{\cA}(\{\prA \mid b\in \prA \})>0$.
	\end{assmpA}
	\Cref{assmp:qSupp}  guarantees that an agent chooses to switch the action with some positive probability $Q(\breve\pi_{a\ast}(\prA))>0$ \textit{whenever }the gross gain $\breve\pi_{a\ast}(\prA)$ is positive. \Cref{assmp:AExAnteavailable} requires any action $b$ to be available with some positive probability; i.e., there exists $\prA$ such that $b\in\prA$ and $\P_\cA(\prA)>0$.  \Cref{assmp:qSupp,assmp:AExAnteavailable} guarantee that a revising agent can always switch to an optimal action with some positive probability; so, a stationary state must be also a Nash equilibrium.
\end{rmk}

\section{Deeper examination behind net gains\label{sec:BehindGain}}
\subsection{Full rationalizability of myopic decision\label{sec:RatMyop}}

The negative correlation between net gains and the dynamics as in property \ref{prop:gzNeg} sounds pretty natural. What is behind it? To answer it, now we relax \Cref{assmp:AIndp_a} and allow $\P_\cA$ to depend on a revising agent's current action $a$. So, given $\x$, let $\P_{\cA a}(\prA;\x)$ be the probability that available strategy set is $\prA$ when a revision agent has been taking action $a$ just before receiving the revision opportunity. This generalization leads us to deeper examination about rational choice in myopic dynamics.\footnote{While our motivation of this generalization here comes from a purely theoretical question, it indeed allows for practically interesting examples as in \Cref{apdx:assumpAInv}.}

Essentially in \textit{evolutionary }dynamics, an agent's choice of a new action is based only on the information about current payoffs $\bpi$. Agents do not take $\x$ into account of the choice; while $\x$ may affect likelihood of sampling of other agents in imitative dynamics or availability of actions in our framework (step 1), it does not affect an agent's decision on a new action (step 2) once the agent finds what actions are available and what payoffs they yield. Furthermore, agents are supposed not to make forecasts about future changes of payoffs. While we could say that an agent might want to maximize a long-term payoff just like in a repeated game, we conventionally justify myopic decision by the notion of a static belief, i.e., assuming that agents somehow believe that the payoff of each action would not change over time. 

To reexamine whether these informational assumptions\footnote{ \citet{Sandholm10Games_PairwiseComparisonDyn}  and \citet[Sec.5.1]{SandholmPopText} clarify what information is needed for an agent to process each of canonical revision protocols.} really justify myopic decision, now we consider an agent who wants to maximize a long-term payoff but only gets information about current state of $\bpi$ and $\x$. Yet, about future payoffs, we assume that agents have static belief, just as in conventional justification. But now we introduce restriction to availability of actions. Then, myopic payoff-based decision becomes not so convincing as a way of rational decision making while we assume limited information. 

The exact optimal action $b_\ast(\bpi)$ may not be available for an agent at a current revision opportunity. For an access to it, an agent may foresee  the next revision opportunity. Then, the agent may choose an action that does not yield the greatest payoff among currently available action but has a better access to the exact optimal action in the next revision opportunity. So, \textit{in general, }an agent should consider not only payoffs $\bpi$ but also accessibility to other actions. Then, our question is under what environment we can still justify the myopic payoff-based decision, so agents would eventually not care about accessibility?

If availability of actions changes with $\x$, then the agent should utilize the information about $\x$ and take it into account of the choice a new strategy, not only $\bpi$. So, \Cref{assmp:AIndp_a}, or a generalized version below for this present set-up, is needed to justify decision based only on payoffs. This assumption implies property \ref{prop:gDx}, i.e., a zero availability effect of $\dot\x$ on $\dot G^\F$.
\begin{assmp}{\ref{assmp:AIndp_x}'}\label{assmp:AIndp_x_gen}
	For any $a\in\cA$, $\P_{\cA a}$ is independent of $\x$.
\end{assmp}

Dependency of availability on an agent's current action suggests that different actions have different accessibility to other actions. However, it matters only if this difference results in different accessibility to \textit{better }actions. To formalize this idea, given payoff vector $\bpi$, let $p^i_{\cA a}$ be the probability that $i$-th best actions are available for an agent with current action $a$; see \eqref{eq:availIth} in \Cref{apdx:Proofs} for a formal definition. A revising agent would not care about the difference in accessibility when deciding on a new strategy, if $p^i_{\cA a}$ does not differ by current action $a$, or more specifically, if the following property holds:\footnote{The difference in lower payoff ranks than both would not be a matter since anyway an agent would not switch to worse actions or a lower rank than the current action, whether it is $a$ or $b$. About ranks between the rank of $a$ and that of $b$---say, $\pi_a<\pi_b$, the access to actions at these intermediate ranks do not make advantage for $a$ over $b$ since these actions are anyway worse than $b$.} 
\begin{description}
	\item[p\label{prop:PIndp_a}] Any two actions $a,b$ have the same accessibility $p^i_{\cA a}=p^i_{\cA b}$ for each payoff rank $i$ higher than ranks of both $a$ and $b$.
\end{description}

This property \ref{prop:PIndp_a} is of not only $\P_{\cA \cdot}$ alone but also of payoff vector $\bpi$, which changes in a game. Then, what assumption on $\P_{\cA \cdot}$ alone can assure property \ref{prop:PIndp_a} for \textit{any} payoff vector $\bpi$, so it is robust to changes in payoffs and holds generally over any games? \Cref{clm:g_pi} tells that it is assured if \textit{and only if }$\P_{\cA \cdot}$ satisfies the assumption below.

\begin{assmp}{\ref{assmp:AIndp_a}'}\label{assmp:AIndp_a_gen}
	For any $a,b\in\cA$ and any $\cA_{ab}\subset \cA\setminus\{a,b\}$ and at any $\x$, $\P_{\cA a}(\{\prA \mid \prA\cap \cA_{ab} \ne\emptyset\};\x)=\P_{\cA b}(\{\prA \mid \prA\cap \cA_{ab}\ne\emptyset\};\x)$.
\end{assmp}

\Cref{assmp:AIndp_a_gen} imposes invariance of availability of actions other than the current action. Specifically, we compare the distribution of available actions $\P_{\cA a}$ when the current action is $a$ and the one $\P_{\cA b}$ when it is $b$. Fix an arbitrary subset $\cA_{ab}'$ of actions other than $a$ and $b$. The assumption requires $\P_{\cA a}$ and $\P_{\cA b}$ to assign the same probability on the event that at least one of the actions in $\cA_{ab}'$ becomes available. So, this is a generalization of \Cref{assmp:AIndp_a}.\footnote{One may ask why we look at a set $\cA_{ab}'$, not each action, and we look at the event that at least one of the actions in the set, not all the actions, is available. A change of the invariance condition in either of these two aspects breaks the equivalence with property \ref{prop:PIndp_a}; see \Cref{apdx:assumpAInv}.}  

Furthermore, \Cref{clm:g_pi} tells that property \ref{prop:PIndp_a} is also equivalent to the following property of $\g_\ast$:
\begin{description}
	\item[g1 \label{prop:g_pi}] Consider arbitrary two actions $a,b\in\cA$; assume $\pi_a\le \pi_b$. Then, i) $g_{b\ast}(\bpi,\x)\le g_{a\ast}(\bpi,\x)$ for any $\bpi,\x$. ii) The inequality is strict if $\pi_a<\pi_b$ and $\0\notin\V_a(\bpi,\x)$.\footnotemark
\end{description}

Property \ref{prop:g_pi} means that a revising agent's switch to a better reply $b$ than the current action $a$ indeed decreases the agent's individual ex-ante net gain from a \textit{next }switch. This is a (stronger) version of property \ref{prop:gzNeg}, which implies decreases in individual ex-ante net gains only on the average of switching agents; so, property \ref{prop:g_pi} is stronger than property \ref{prop:gzNeg}. \Cref{clm:g_pi} suggests that, if \Cref{assmp:AIndp_a_gen} or property \ref{prop:PIndp_a} does not hold, a switch to a better reply based on $\bpi$ may not decrease the remaining gain; so, agents who take the best available action may not achieve the greatest exploitation of the remaining individual ex-ante net gains. On the other hand, if \Cref{assmp:AIndp_a_gen} is satisfied, this discrepancy disappears. 

\begin{thm}\label{clm:g_pi}
	A cost-benefit rationalizable dynamic $\V$ satisfies the following. 
	\begin{enumerate}[a)]
		\item The following statements are equivalent to each other: 
		\begin{itemize}
			\item $\P_{\cA\cdot }$ satisfies \Cref{assmp:AIndp_a_gen}. 
			\item Property \ref{prop:PIndp_a} holds for any $\bpi,\x$. 
			\item Property \ref{prop:g_pi} holds for any $\bpi,\x$. 
		\end{itemize}
		\item Property \ref{prop:g_pi} implies property \ref{prop:gzNeg}.
	\end{enumerate}
\end{thm}
\footnotetext{The latter condition $\0\notin\V_a(\bpi,\x)$ can be removed under \Cref{assmp:AExAnteavailable_gen,assmp:qSupp}, since $\pi_a<\pi_b$ implies $a\notin b_\ast(\bpi)$ and then these assumptions imply $\0\notin\V_a(\bpi,\x)$.}

In sum, \Cref{assmp:AIndp_a_gen,assmp:AIndp_x_gen} for our notion of ``full rationalizability" indeed rationalize myopic payoff-based decision. By the way, with generalization of \Cref{assmp:AExAnteavailable} as below, all these generalized assumptions replace all the preceding assumptions.
\begin{assmp}{\ref{assmp:AExAnteavailable}'}\label{assmp:AExAnteavailable_gen}
	For any $a\in\cA$ and any $b\in\cA\setminus\{a\}$, $\P_{\cA a}(\{\prA \mid b\in \prA \})>0$.
\end{assmp}

\begin{prop}\label{clm:gen}
	By replacing Assumptions \ref{assmp:AIndp_a}, \ref{assmp:AIndp_x}, \ref{assmp:AExAnteavailable} with  Assumptions \ref{assmp:AIndp_a_gen}, \ref{assmp:AIndp_x_gen}, and \ref{assmp:AExAnteavailable_gen} respectively, all the theorems so far still hold.
\end{prop}

\subsection{Cause of instability in imitative dynamics}\label{sec:ImitDyn}
The replicator dynamic can be constructed as imitative dynamics \citep{Schlag98JET_Imitate}. In this dynamic, only Lyapunov stability (no farther escape) is obtained for stable games in general. For asymptotic stability (convergence in addition to no father escape), \textit{strict }static stability of a game is needed. Requirement of strictness of stable games might  appear to be a subtle, minor difference. However, we should notice that it excludes a game with a constant payoff vector $\F(\x)\equiv\bpi^0$ at any $\x\in\Delta^\cA$. In such a situation, we naturally expect economic agents to switch to the best actions, which do not change over time in this situation; the population state should converge to a Nash equilibrium, as long as agents' choices are driven by incentives. This is assured for fully rationalizable dynamics. In contrast, imitative dynamics result in a cyclic trajectory, staying away from the equilibrium. This is generally proven for null stable games, which include zero-sum games, where \Cref{eq:DF_StatStblCondn} always holds with equality \citep{AkinLosert84JMathB_EvolDyn_ZeroSum,Hofbauer96JMathB_EvolDyn_Bimatrix_Hamiltonian}; see \citet[Sec.~9.1]{SandholmPopText}.

What is behind this non-stability of imitative dynamics? The net gain gives an answer. Imitative dynamics can be fit into our framework and thus cost-benefit rationalizable by defining the probability that action $b$ is available, i.e., $\P_{\cA a}(\{\prA \mid b\in \prA\})$ as the probability that a $b$-player is sampled from the society. The essence of such imitative dynamics is that the sampling probability depends on the current actual distribution of actions in the society. Hence, Assumption \ref{assmp:AIndp_x} does not hold and thus imitative dynamics are not fully rationalizable. While individual agents switch to better replies, these switches may not decrease the remaining net gains. 

Actually, imitation refuels net gains. Imagine a deviation from an equilibrium and look at actions that have greater shares of players than in the equilibrium. While static stability yields a negative payoff effect of $\dot\x$ on $\dot G^\F$ through lowering relative payoffs of these actions, imitation rather generates a positive availability effect as these strategies become more likely to be sampled and thus more available. This refuels the ex-ante net gain. This intuition is supported by the following theorem.

\begin{thm}\label{clm:ReplRefuelG}
	Let $\V^\F$ be a  replicator dynamic in a zero-sum game. Then, the availability effect satisfies the following:
	$$ \x\cdot \pdif{\g_\ast}{\x} (\x,\bpi) \dot\x \begin{cases}
		\ge 0 &\text{(i) for any $\dot\x\in\V^\F(\x)$ at any $\x\in\Delta^\cA$},\\
		>0  &\text{ (ii) if $\x\notin\RP(\V^\F)$}.
	\end{cases}$$
\end{thm}

\section{Extensions of the net gain approach\label{sec:Ext}}

\subsection{Rationalizability and gains in other major dynamics}\label{sec:Modified}
The ``generality'' of our approach lies in the idea itself---to be summarized in \Cref{sec:Concl}; the analysis based on the ``basic framework' is made to convince readers of mathematical validity of our approach and to present a concrete way to construct the net gain function and prove crucial properties.

The essence of our approach is that the aggregate ex-ante net gain generally works as a Lyapunov function to prove or disprove dynamic stability as long as we can reconstruct a dynamic in a cost-benefit rationalizable way, i.e., based on optimization by making explicit costs and constraints that prevent agents from exact optimization. This idea is applicable to \textit{perturbed best response dynamics} (pBRD), where agents follow a kind of quantal response rules such as logit choice rules. According to the optimization-based construction of pBRD by \cite{HofbauerSandholm02EMA_StochFictPlay}, there are two ways to make it cost-benefit rationalizable. First, each agent is considered as facing idiosyncratic payoff perturbations while the agent chooses the exact optimal choice given the \textit{perturbed }payoff vector. The perturbed BRD emerges as the aggregation of such heterogeneous exact BRD. We discuss heterogeneous dynamics in \Cref{sec:hetero} as an extension of the basic model to games with finitely many distinctive populations. Second, an agent chooses mixed strategy $\y\in\Delta^\cA$ and there is a control cost $V(\y)$ to prevent the agent from taking a pure strategy. The control cost perturbs the agent's optimal strategy. We can interpret the known Lyapunov function for perturbed BRDs in \citet[Theorem 3.1]{HofSand07JET_RdmnDistPayoff} as the aggregate net gain function after counting the control cost.

Cost-benefit rationalizability of \textit{excess payoff dynamics }needs a little modification of our framework. In those dynamics, an agent's switching rate to a new action is an increasing function of its excess payoff, i.e., the relative payoff of the action compared with the \textit{average payoff in the population} $\bpi\cdot\x$. To rationalize it, we regard the population's average strategy, i.e., $\x$, as the mixed-strategy status quo for a revising agent: if an agent chooses the status quo to avoid the payment of a switching cost, the new action is chosen randomly according to the current action distribution $\x\in\Delta^\cA$ in the population; a revising agent can be interpreted as a newborn agent, replacing an old agent. Then, the excess payoff determines the switching rate. See \Cref{apdx:Modified} for details and the construction of the net gain function. 

Evolutionary dynamics may be built upon complex decision processes. For example, there has been a recent development in ``sampling dynamics"\footnote{See \citet{OSTsBRD,Sandholm_Izquierdo_19TE_Stabl_BRExpDyn,Sandholm_Izquierdo_19JET_Stabl_BRExpDyn,Arigapudi_etal_21_JET_InstblPD_BEPD,Sawa_Wu_21_StatInfer_EvolDyn}, for recent examples.}. Basically they consider the BRD without direct observation of current payoff vector $\bpi$, but with inference of $\bpi$ from sampling of other agents' earned payoffs or actual choices of actions. This sampling-inference process can be considered as a conversion of $\F$. Let $\textbf S:\bpi\mapsto\hat\x$ be the action sampling process that generates the inferred action distribution $\hat\x$ from observations of actions given the true action distribution $\x$. Then, it converts the game $\F$ to $\hat\F(\x)=(\F(\textbf S (\x))$.  Let $\textbf S:(\bpi,\x)\mapsto\hat\bpi$ be the payoff sampling process that generates the inferred payoff vector $\breve\bpi$ from observations of payoffs given the current payoff vector $\bpi$ and social state $\x$. Then, it converts the game $\F$ to $\hat\F(\x)=\textbf S(\F(\x),\x)$. Then, the approach here is applied to these sampling-based dynamics by regarding them as the standard BRD in such a \textit{converted} game $\hat\F$.

\subsection{Heterogeneity}\label{sec:hetero}
The model and all the propositions can be extended to a multi-population game in which (finitely many) different populations may have different payoff functions and/or different revision protocols, as long as each population's revision protocol is cost-benefit rationalizable.\footnote{\cite{Zusai_HeteroEvolDyn} allows the number of types to be continuously many and provides the measure-theoretic definition of heterogeneous evolutionary dynamic, while clarifying regularity conditions for the unique existence of a solution trajectory and proving dynamic stability of equilibria in potential games.} 
We obtain the Lyapunov function for a multi-population game just by adding up agents' ex-ante net gains over all the populations.\footnote{
	It can be cost-benefit rationalizable in the basic framework in this paper or in some modified ways as argued in \Cref{sec:Modified}. In particular, we can consider a hybrid dynamic where an agent may keep the current choice of an action as a status quo or may be replaced with a new agent. Further, even if a population's dynamic is not cost-benefit rationalizable, we can include it as long as we can somehow find functions $G$ and $H$ that satisfy the properties in \Cref{clm:GZero,clm:DG} for this dynamic.} 
It is easy to see properties of $G$ are preserved by summation; see \Cref{apdx:hetero}.  

The extension to a multi-population setting covers some interesting class of games. One is a \textit{saddle game}, which is a generalization of zero-sum games and also related to the literature on robust equilibria under incomplete information by \cite{NoraUno14JET_SaddleFn_RobustEqia}.\footnote{Evolutionary dynamics in this class of games can date back to \cite{Kose56EMA_SolSaddle_DiffEqm} that considers computation of a saddle point problem through evolutionary algorithms.} Another is an \textit{anonymous game}, in which agents' payoffs depend only on the aggregate action distribution over the society.\footnote{Aggregative games \citep{Jensen18Hdbk_AggregativeGames} such as Cournot competition are anonymous games, though the former assumes that strategy variable $a\in\cA$ represents some numerical value and it makes sense to add them over all agents as $\sum_{a\in\cA} a x_a$. A normal-form game in uniform random matching is an anonymous game but not necessarily an aggregative game.} The aggregation represents unobservability of each agent's affiliating population when they play the game and their payoffs are assessed. 

\cite{SandholmCongstPr,SandholmNegImpl} proposes to extend the idea of Pigouvian tax to the evolutionary dynamic setting in order to achieve the maximal aggregate payoff at the limit state under negative externality, with the central planner knowing only the aggregate action distribution. To prove global convergence to the social optimum, the aggregate payoff functions is supposed to be concave; the Pigouvian tax makes it a potential function in a game with tax and thus concavity of the potential function implies global convergence. Since a stable game is generalization of a concave potential game and removes the assumption of exact symmetry of externalities, our result suggests robustness of Sandholm's idea of ``evolutionary implementation" of the social optimum.

\subsection{Discrete-time finite-agent models}\label{sec:discrete}
Adaptive dynamics may be modeled in a discrete time horizon with a population of finitely many agents. It may be a popular way to model dynamics in the literature of learning dynamics \citep{Young04} and also a way to enable numerical computation as in agent-based simulations \citep{Izquierdo13JStatPhys_StochApprox_SimpleSimulModels,Izquierdo_Sandholm_AgentBasedEGD} or to fit into practical applications  as in distributional controls \citep{MardenShamma15hdbk,Quijano17IeeeCtrlSys_RolePopGameEvolDyn_DistCtrlSys}.\footnote{On stochastic adaptive learning in finite normal-form games, \cite{Funai_18ET_Conv_StochAdaptDyn} proposes a framework to unify various learning models such as stochastic fictitious play, experience-weighted attraction learning, etc.; then, he investigates how convergence to a quantal response equilibrium depends on these learning parameters in several classes of games. Yet, he does not use a Lyapunov function approach.} 
 A large population model over a continuous time horizon is considered as a limit case of such discrete models. To link them, our approach has two remarkable features:	 
\begin{enumerate}[a)]
	\item We construct a differential inclusion $\dot\x\in\V^\F(\x)$ from optimization. Because of this, the correspondence $\V^\F:\Delta^\cA\rightrightarrows\R^\cA_0$ is non-empty, convex-valued and upper semicontinuous.
	\item The Lyapunov function is used to prove equilibrium stability.
\end{enumerate} 
These features enable us to utilize findings in the literature on discrete dynamics:
\begin{enumerate}[ a)]
	\item According to \citet{RothSandholm_StochApproxDI}, the differential inclusion indeed \textit{approximates} the middle-run behavior and the limit states of a finite-population dynamic, possibly over a discrete time horizon.
	\item According to \cite{EllisonFudenImhof16JET_FastConvEvol_Lyapunov}, the Lyapunov function guarantees that the expected convergence time is bounded.
\end{enumerate} 
The first property assures a close linkage between our continuous dynamics and computational discrete dynamics; simulations serve a legitimate approximate substitute for mathematical analysis and mathematical analysis tells fundamental properties of simulation models. The second property relieve computational researchers' worry by assuring that computation should finish in a bounded time (with an arbitrary small error bound).  

\section{Discussion: why not other ways?}\label{sec:dis}

\subsection{Why not other Lyapunov functions?}

\subsubsection*{The gross gain}\label{sec:GrossGain}
To calculate a net gain, we need switching cost $q$. One may think it artificial. Why can't we simply use gross gains? The first and obvious reason is that a zero gross gain may not imply a rest point of a dynamic. It is the case when there is a positive lower bound  of switching costs so $Q(\bar q)=0$ with some $\bar q>0$.\footnote{Deterministic (non-zero) switching costs are considered in  \cite{Norman09GEB_RapidEvol,Norman10ThDc_CycleEqmEvol} on stochastic evolution and in  \cite{LipmanWang00JET_SwitchCostFrequentRepeat,LipmanWang09GEB_SwitchInfRepeat} on repeated games.} Then, \Cref{assmp:qSupp} is violated and there are non-equilibrium rest points, i.e, $\NE(\F) \subsetneq \RP(\V^\F)$. The aggregate gross gain is defined as\footnote{Rigorously speaking, Lipschitz continuity is needed to use $\Gamma^\F$ as a Lyapunov function; see \Cref{thm:Lyapunov_DI} in \Cref{apdx:DynMath}. It is retained by replacing $Q$ with $Q_-$. Rigorous and extensive study of $\Gamma^\F$ was done in an older working paper version (Sec. 6.1) of the present paper, available at \url{https://arxiv.org/pdf/1805.04898v4.pdf}.}  
$$\Gamma(\x,\bpi)\coloneqq \sum_{a\in\cA} x_a \sum_{\prA\subset\cA} \P_\cA (\prA;\x) Q(\breve\pi_{a\ast}(\prA)).$$ 
This attains zero if and only if $\x\in\Delta^\cA(b_\ast(\bpi))$; thus, given game $\F$, $\Gamma^\F$ attains zero at Nash equilibria but not a non-equilibrium rest points. Yet, the net gain works since it does not need \Cref{assmp:qSupp} to have $\RP(\V^\F)={G^\F}^{-1}(0)$ as in property \Cref{prop:GZero}. 

Second, even if a dynamic satisfies \Cref{assmp:qSupp}, the gross gain does not provide a close link between static and dynamic stability unlike the net gain. Compare the following with property \ref{prop:GDpi}.
\begin{align*}
	\pdif{\Gamma}{\pi}(\x,\bpi) \Delta\bpi 
	&= \sum_{a\in\cA} x_a \sum_{\prA\subset\cA}\P_{\cA a}(\prA)  Q(\breve\pi_{a\ast}(\prA))\Delta\breve\pi_{a\ast}(\prA)\\
	&\qquad +\sum_{a\in\cA} x_a \sum_{\prA\subset\cA}\P_{\cA a}(\prA)  \{Q'(\breve\pi_{a\ast}(\prA))\Delta\breve\pi_{a\ast}(\prA)\} \breve\pi_{a\ast}(\prA),
\end{align*}
The first term is equal to $\Delta\x\cdot\Delta\bpi$ with $\Delta\x\in\V(\x,\bpi)$, same as property \ref{prop:GDpi}. To understand the second term, observe that a change in payoff vector by $\Delta\bpi$. causes a change in the conditional switching rate $Q(\breve\pi_{a\ast}(\prA))$ by $Q'(\breve\pi_{a\ast}(\prA))\Delta\breve\pi_{a\ast}(\prA)$. This reflects the changes in switching decisions for agents who are on the margin between choosing to switch and not. By changing from no switch to switch, each one of such marginal agents adds the individual gross gain $\breve\pi_{a\ast}(\prA)$ to the aggregate gross gain $\Gamma$. Notice that $\breve\pi_{a\ast}(\prA)$ may be significantly (not marginally) positive even for those marginal agents, since the \textit{gross} gain is not an exact determinant for switching decision unlike the net gain. Because of this, the gross gain may not decrease even under a fully rationalizable dynamic in a strictly stable game; see \Cref{exmp:RPS} in \Cref{apdx:dis}. Thus, the direction of the change in $\Gamma^\F$ depends on specifications of dynamics such as $Q'$ beyond cost-benefit rationalizability.

\subsubsection*{A potential of a game}
Now let's compare the net gain with other existing Lyapunov functions in the preceding literature. \citet{Sandholm01JET_Potential} imports the notion of a potential game \citep{MondererShapley96GEB_Potential} to a population-game setting and proves that a potential function $f$ of a potential game works generally as a Lyapunov function for a wide class of evolutionary dynamics. This approach indeed requires quite a small assumption on $\V$, only a positive correlation between the relative payoff of each action and the change in the mass of players of the action, i.e., $\F(\x)\cdot\dot\x\ge 0$. Since $\dot f=\F(\x)\cdot\dot\x$, this positive correlation assures that the potential increases over time until it reaches a local maximum, which is known to be a Nash equilibrium. But the existence of a potential function is not taken as granted for any game. Under continuous differentiability of $\F$, its existence (the definition of a potential game) is equivalent to exact symmetry of externality between two actions, i.e., $\tpdif{F_a}{x_b} \equiv \tpdif{F_b}{x_a}$ for any pair of actions $a,b$. Since it is an equality condition (unlike the definition of static stability), small perturbation in payoffs may break this condition. 

Notice that the correlation $\F(\x)\cdot\dot\x$ can be interpreted as the aggregate change in agents' payoffs, i.e., the aggregate gross gain. The first and second order condition for a local maximum $Df(\x^\ast)=0$ and $D^2f(\x^\ast)$ being negative semidefinite imply that monotone increase in $f$ to a local maximum is equivalent to monotone decrease in $\Gamma^\F$ to zero. Thus, the gross gain may work as a substitute for a potential function and may be applicable beyond potential games, though applicability (monotone decrease) depends on more details of $\V$ than the net gain as argued above.

\subsubsection*{Geometric distance or divergence}
Another popular candidate for a Lyapunov function is a geometric distance from an equilibrium.  \citet{MerikopoulousSandholm18JET_RiemannianGameDyn} show that a generalized distance based on Riemannian geometry works to prove asymptotic stability of the unique Nash equilibrium in a strict stable game (implying the uniqueness) under a wide range of imitative dynamics. The known Lyapunov function for the replicator dynamic, the Kullback-Leibler divergence is such a distance (divergence) based on a Riemannian metric.\footnote{The KL divergence has become also popular in the recent literature on misspecified learning under incomplete information to characterize an equilibrium; for example,  \citet{Fudenberg_Lanzani_Strack_21Ecta_LimitPts_EndogMisspecifLearn,Massari_Newton_20_LearnEqm_Misspecif}. The true objective probability distribution is naturally taken as the origin to measure the KL divergence.}

This geometry-based approach is an advantage of this geometric approach that it works to prove equilibrium stability under imitative dynamics in a strictly stable game, compared to our net gain approach. However, it has a crucial limitation in practical use. That is, to define a distance-based Lyapunov function, we first need to identify a limit state as the origin to measure the distance. The strictness of a stable game implies the uniqueness of a Nash equilibrium and thus gets rid of this limitation. 

The net gain approach does not require prior knowledge about a limit state, or even whether the dynamic has a limit state; it works to explain instability and also can be applied to a set of rest points. It would be practically helpful in computational study of new dynamics or complicated games, since a practitioner may not know even where an equilibrium is and the practitioner may want to use evolutionary dynamics (or adaptive agent-based dynamics) in order to numerically identify equilibria.\footnote{For example,  \citet{Lamotte_Geroliminis_21TranspResB_Monot_TripSch} use the net gain to numerically confirm convergence in trip scheduling problems in transportation engineering, citing a working paper version of the present paper. }

\subsection{Why not assuming payoff monotonicity for stability?}\label{sec:OrderComp}
Full rationalizability would appeal to economic intuition and also guarantees the link between static and dynamic stability without needing additional assumptions; both \Cref{assmp:AIndp_a_gen,assmp:AIndp_x_gen} can be understood as a basis to make a myopic decision based on $\bpi$ compatible with true incentives for revision as measured by $\g_\ast$. This marks a stark difference from the preceding attempts to complete this link in evolutionary game theory. In this field, \textit{payoff monotonicity} has been regarded as a notion of incentive compatibility of a dynamic. It means that a strategy that yields a greater payoff than another strategy should attract more agents (either in the number/mass or in the rate of increase) than the other. \citet[Counterexample 2]{Friedman91Ecta_EvolGamesEcon} finds a counterexample that payoff monotonicity does not guarantee convergence to an ESS; a similar counterexample is also reported by H\&S (Example 6.1).\footnote{In \Cref{apdx:dis}, we reexamine Friedman's example and find that his dynamic is indeed not fully rationalizable. We also find that a fully rationalizable dynamic may not be payoff monotonic.} More broadly, the replicator dynamic is payoff monotone in the growth rate but fails to guarantee convergence to an equilibrium under (null) stable games, as argued in \Cref{sec:ImitDyn}.

The preceding literature tried to fix it by imposing an \textit{additional }assumption. For excess payoff dynamics, H\& S proposed \textit{integrability }of a function that maps an excess payoff vector to the rate of switching to each action.\footnote{The idea stems from the literature on Hannan consistency in repeated games: \citep{FudenLevine98_Learinng,Hart_MasColell_01JET_GenClass_AdaptiveStr}.} Integrability means the existence of a potential of this switching rate function just like a potential of a game.\footnote{Integrable dynamics may not be fully rationalizable, since the former allows agents to switch to actions that yield lower payoffs than the average. Full rationalizability does not imply integrability. See \cref{apdx:Ext}.} However, there is no single principle to cover both the payoff monotonicity condition and the integrability condition. Anyway, the integrability condition is only for excess payoff dynamics. Besides, it is hard to give an economic meaning of integrability and economic intuition behind why it works.

\citet{Sandholm14DGA_Integrability_GameDyn} admits this point\footnote{\cite[p.96]{Sandholm14DGA_Integrability_GameDyn} says ``However, despite the usefulness of integrability for proving convergence results, the game-theoretic interpretation of the integrability of choice functions is obscure."} and tries to find an economic meaning of the integrability assumption. For this, he focuses on a \textit{single-player} decision model where the payoff vector changes \textit{randomly} over time, and introduces ``signals'' to tell \textit{future }performances; then, integrability means \textit{no }correlation between switching decisions and those signals. So, switching decision is truly myopic. This characterization of integrability may accord with our characterization of full rationalizability, especially \Cref{assmp:AIndp_a_gen}. In \Cref{clm:g_pi}, we find this assumption is equivalent to property \ref{prop:g_pi}; this property implies that the ex-ante net gains $\g_\ast$, which account for accessibility to other actions, do not carry information beyond payoff vector $\bpi$ for an agent's decision on a new strategy.%

\section{Concluding remarks}\label{sec:Concl}
In this paper, we present a unifying and economically intuitive approach to link static and dynamic stability based on cost-benefit rationalizability behind agents' updates of choices. The approach is summarized as follows:
\begin{enumerate}
	\item Rationalize an individual agent's decision rule (revision protocol) based on optimization; distortions from exact optimization should be rationalized with costs and constraints in strategy revision. 
	\item Then, define a \textit{net} gain as the payoff improvement after deducting those additional costs. Calculate the ex-ante net gain as the greatest available net gain subject to those constraints. Aggregate it over all the agents. It works as a Lyapunov function.
\end{enumerate}
Our general approach not only provides a general proof of dynamic stability with mathematically universal and economically intuitive logic but also explains causes of instability and extends it to the heterogeneous setting or to more complex dynamics. Compared to the other known Lyapunov functions, the net gain has two distinctive advantages. First, it works even if rest points of a dynamic do not coincide with equilibria of a game. Second, calculation of the net gain does not require prior knowledge about rest points or a limit state.

Our notion of cost-benefit rationalizability is new in evolutionary game theory, but prevalent in any other areas of economics. The optimization principle of ``rational choice'' is so powerful to reveal the underlying preference, which in turn allows us to test rationality empirically (e.g.  \citealt{Kitamura_Stoye_18Ecta_NonparaAnal_RandUtilModels}). In estimation of discrete choice models, econometricians attempt to explain empirical data with random utility models \citep{AndersonDePalmaThisse92_DiscChoice}. In game experiments, quantal response equilibrium \citep{McKelveyPalfrey95JET_QRE} has been dominantly used as a solution concept; the quantal response is regarded as optimization under idiosyncratic payoff perturbation. In recent development of decision theory, theorists try to rationally axiomatize behavioral distortions (e.g. \citealt{Natenzon_19JPE_RndmChoice_Learn}).
Those economists admit distortions from exact optimization but try to understand each of those distortions by revealing hidden costs, constraints and structures in decision making, and have proven usefulness of the rational-choice approach by deriving rich implications and presenting good matches with empirical data. Our approach shares this spirit and proves usefulness in dynamic modeling by solving a fundamental problem of finding a link between static and dynamic stability.

While our analysis is confined here to a normal-form population game, the author believes that there is a broad range of future applications of our approach  to proving equilibrium stability beyond this simple setting; see  \cite{SawaZusai_multitaskBRD} for an example of our approach to the BRD with a nested choice structure in multitasking environments. Typically in sequential-move, repeated or coalitional games, one may define a version of static stability that is tailored to each specific situation in order to refine an equilibrium concept.\footnote{See \citet{GarciaVeelen16JET_InOutEqm1_RepeatDisc,GreveOkuno09RES_VoluntSeparableRepeatedPD,GreveOkuno16_DiversePatterns_VoluntaryPartner} for examples in repeated games. Models of preference evolution  normally take the ``indirect evolution approach'': see  \citet{DekelElyY07RES_EvolPref,AlgerWeibull13EMA_HomoMoralis}. See \cite{Demuynck_etal_19Ecta_MyopicStblSet_SocialEnvir} for an example of a general definition of static coalitional stability based on myopic adjustments.} Equilibrium refinement with those versions of static stability imposes a condition on the payoff function (the game) alone, while the adjustment process may be discussed only casually. One may hope to justify the ad-hoc refinement by establishing dynamic stability. This paper suggests to formally define such a dynamic from the casual idea about the underlying adjustment process and track down changes in the remaining \textit{net gain} under the dynamic; then dynamic stability should be obtained.

Finally, it is notable that our net gain approach to analysis and prediction of dynamic stability is based on micro-foundation approach to construction of evolutionary dynamics, as manifested in \cite{SandholmPopText} and materialized in his numerous papers. The book presents an organized, rigorous and universal procedure to define an evolutionary dynamic from formulation of individual decision rules (revision protocols). Then, it would be natural as \textit{a next step }to explain these revision protocols \textit{economically}, i.e., rationalize various styles of bounded rationality by costs and constraints in revisions. This rationalization is what we push further in the current paper; and it leads to the idea of net gains as a micro-founded measuremment of convergence and also the notion of full rationalizability of myopia as a micro-level condition for establishing the link between static and dynamic stability. In this paper, we succeed the spirit of  \cite{SandholmPopText} and find a ground on the field of economics to plant the spirit and to harvest its fruits. Therefore, this paper is truly dedicated to Bill Sandholm, not solely for our long friendship but for his clear vision of EGT.

\appendix
\setcounter{section}{0}
\renewcommand{\thesection}{A\arabic{section}}
\section{Appendix to \Cref{sec:Prelim} \label{apdx:Prelim}}
\subsection{Brief summary of math on dynamic systems \label{apdx:DynMath}}
See \cite{Smirnov01DI} for a general extensive reference on diffierential inclusions.  For a rigorous but handy quick reference, see \cite[Ch. 7 Appendix]{SandholmPopText}. 

\paragraph{Concept of a solution trajectory.}
Since a differential inclusion may allow for discontinuous switches of transition vectors, a solution trajectory may exhibit non-smoothness. While there are several concepts of solution trajectories, we adopt a Carath\'{e}odory solution: $\{\x_t\}_{t\ge 0}$ is a Carath\'{e}odory solution trajectory for $\dot\x\in\V^\F$ if it is Lipschitz continuous in $t$ at every $t\ge 0$ and also differentiable with derivative $\dot\x_t\in \V^\F(\x_t)$ at almost every $t$. The existence of a Carath\'{e}odory solution trajectory from an arbitrary initial point $\x_0\in\Delta^\cA$ is assured under our framework, because $\V^\F$ of a cost-benefit rationalizable dynamic (either in the basic set-up in \Cref{sec:Basic} or in the modified one in \Cref{apdx:Modified}) is a bounded and upper semicontinuous correspondence with compact and convex values \citep[Ch.~4]{Smirnov01DI}. Unlike a differential equation, this regularity condition does not necessarily mean uniqueness; multiple solution trajectories are not uncommon.

\paragraph{Stability concepts.} \textbf{Lyapunov stability} of a set $X^\ast$ (possibly a singleton of a point) means that, even if the state deviates from $X^\ast$, the state never goes too far away from $X^\ast$ under any solution trajectory. If the state returns to $X^\ast$ asymptotically, we say $X^\ast$ is \textbf{attracting}. It is \textbf{asymptotically stable} if it is both Lyapunov stable and attracting.

\paragraph{Lyapunov function.} 
In this paper, we refer to the following version of the Lyapunov stabiliity theorem.\footnote{See \citet[Theorem 8.2]{Smirnov01DI} for a standard version that applies to stability of a single state. \citet[Theorem 8.2]{ZusaiTBRD} modifies it for stability of a set. \Cref{thm:Lyapunov_DI} adopts the former for local stability of an isolated rest point and the latter for global stability of a set.}
		
\begin{thm}[Lyapunov stability theorem]\label{thm:Lyapunov_DI}
Consider a differential inclusion $\V^\F$ on $\Delta^\cA$. Suppose that continuous function $W:\Delta^\cA \to\R$ and lower semicontinuous function $\tilde W:\Delta^\cA \to\R$ satisfy (a) $W(\x)\ge 0$ and $\tilde W(\x)\le 0$ at all $\x$. 
	\begin{description}
	\item[Global stability.] Let $X^\ast$ be a non-empty closed set in $\Delta^\cA$. Assume that (b) $W^{-1}(0)= \tilde W^{-1}(0)=X^\ast$ and (c) $W$ is Lipschitz continuous. If 
		\begin{equation}
			\dif{}{\x}W(\x)\dot\x \le \tilde W(\x) \qquad \text{for any }\dot\x\in \V^\F(\x) \label{eq:LCondtn_DI}
		\end{equation}
	whenever $W$ is differentiable at $\x$, then $X^\ast$ is globally asymptotically stable under $\V^\F$.
	\item[Local stability.] Let $\x^\ast$ be the only rest point in a neighborhood $X'$. Assume that  (b') $\x=\x^\ast \Leftrightarrow  [W(\x)=0  \text{ and } \x\in X'] \Leftrightarrow  [\tilde W(\x)=0  \text{ and } \x\in X']$ and (c') $W$ is Lipschitz continuous in $X'$. If \eqref{eq:LCondtn_DI} holds whenever $W$ is differentiable at $\x\in X'$, then $\x^\ast$ is asymptotically stable under $\V^\F$.
	\end{description}
We call $W$ a \textit{Lyapunov function} and $\tilde W$ a \textit{decaying rate function}.
\end{thm}
	
\section{Proofs of theorems\label{apdx:Proofs}}
\subsection{Proof of \Cref{clm:GZero,clm:DG,clm:g_pi}}
While we organize theorems in the main text from the most appealing ones to deeper ones, here we drastically reorganize the proofs of theorems in a logical order, i.e., starting from basic properties of individual net gain function $\g_\ast$ and then deriving corresponding properties of aggregate net gain function $G$. In the proofs below, we adopt a general set-up in which availability of actions $\P_{\cA \cdot}$ may depend on current action of a revising agent. So, these proofs reflect \Cref{clm:gen}. \footnote{In these proofs, we omit $\x$ from arguments of $\P_{\cA \cdot}$ since $\x$ is fixed in all the statements except property \ref{prop:gDx} (while $\Delta\x$ from a fixed $\x$ may be involved) and, for this property, we assume \Cref{assmp:AIndp_x_gen}.}

\subsubsection*{Basic properties of $\g_\ast$}
Property \ref{prop:GZero} in \Cref{clm:GZero} and property \ref{prop:GDpi} in \Cref{clm:DG} hold for the aggregate net gain function $G$ for any cost-benefit rationalizable dynamics. They come from the following properties of the individual ex-ante net gain function $\g_\ast$.

\begin{lem}\label{clm:g0}
	Consider a cost-benefit rationalizable dynamic $\V$. Then, the individual ex-ante net gain function $\g_\ast:\Delta^\cA\times\R^\cA \to\R^\cA$ satisfies the following properties. 
	\begin{description}
		\item[g0-a) \label{prop:gZero}] i) $g_{a\ast}(\x,\bpi)\ge 0$ for any $\x\in\Delta^\cA,\bpi\in\R^\cA;$ ii) $g_{a\ast}(\x,\bpi)=0$ if and only if $ \0\in\V_a(\x,\bpi) $.
		\item[g0-b) \label{prop:gDpi}] $\g_\ast$ is Lipschitz continuous and  thus differentiable almost everywhere. Whenever  $g$ is differentiable at $(\x,\bpi)\in\Delta^\cA\times\R^\cA,$ we have
			\begin{equation}
			\pdif{g_{a\ast}}{\bpi}(\x,\bpi) \Delta\bpi=\z_a\cdot\Delta\bpi \quad\text{ for any }\Delta\bpi\in\R^\cA, \z_a\in\V_a(\x,\bpi). \label{eq:Dg_pi}
		\end{equation}
	\end{description}
\end{lem}

\begin{proof}[Proof of \Cref{clm:g0}]
\noindent\textbf{\ref{prop:gZero} } Part i) is immediate from the definition of $g_{a\ast}$. Note that the term $\E_Q[\breve\pi_{a\ast}(\prA)-q]_+$ reflects the fact that a revising agent switches to a different new action only if $\pi_\ast(\prA)\ge \pi_a+q;$ so does this non-negativity of $g_{a\ast}$.
	
ii) By the definition of $\V_a$, the stationarity condition $\0\in\V_a(\bpi,\x)$ is equivalent to $0\in\cQ(\breve\pi_a(\prA))$  whenever $\P_{\cA a}(\prA)>0$.\footnote{There is a possibility that the stationary condition is satisfied by $\e_a\in \Delta^\cA(b_\ast(\bpi,\prA))$, which is equivalent to $a\in b_\ast(\bpi,\prA)$. This implies $\breve\pi_a(\prA)=0$. Since $q$ takes only a non-negative value, it is followed by $Q(\breve\pi_a(\prA))=0$.} This statement $0\in\cQ(\breve\pi_a(\prA))$ means $Q_-(\breve\pi_a(\prA))=\P_Q(q< \breve\pi_a(\prA))=0$ and thus it is equivalent to $\E_Q[\breve\pi_a(\prA)-q]_+=0$. Hence, the stationarity condition is equiavelent to $\E_Q[\breve\pi_a(\prA)-q]_+=0$ whenever $\P_{\cA a}(\prA)>0$, which holds if and only if $g_{a\ast}(\bpi,\x)=0$ by the definition of $g_{a\ast}$.

\noindent\textbf{\ref{prop:gDpi} } First of all, let $\fb(\prA)\in \prA$ be a choice of a (not necessarily optimal) action from non-empty available action set $\prA\subset\cA$; we could call $\fb:2^\cA \setminus\emptyset\to\cA$ a policy function for a revising agent to make a choice after observing (non-empty) available action set. Let $\fB$ be the set of all such policy functions; notice that it is a finite set as long as $\cA$ is a finite set. Define the ex-ante gain from policy $\fb$, $g_{a\fb}$, by
	$$ g_{a\fb}(\x,\bpi) \coloneqq \sum_{\prA\subset \cA} \P_{\cA a}(\prA) \E_Q [\pi_{\fb(\prA)}-\pi_a-q]_+.$$
	If an agent follows the randomly constrained optimization protocol \eqref{eq:OptProtocol}, $\fb(\prA)$ should be chosen from $b_\ast(\prA)$ in each possible realization of $\prA$; thus,
	$$ g_{a\ast}(\x,\bpi) = \max_{\fb\in\fB} g_{a\fb}(\x,\bpi).$$
	Specifically, the maximum is attained by a policy $\fb_\ast\in\fB$ such that $\fb_\ast(\prA)\in b_\ast(\prA)\subset \prA$.
	
	With $\fb$ fixed arbitrarily in $\fB$, function $g_{a\fb}$ is Lipschitz continuous and thus differentiable almost everywhere. It is differentiable at $\bpi$ if and only if $\P_Q$ is continuous at $\pi_{\fb(\prA)}-\pi_a$ of all $a\in\cA$ and all $\prA\subset\cA$ with $\P_{\cA a}(\prA)>0$. Then, the derivative is 
	$$ \pdif{g_{a\fb}}{\pi_b}(\x,\bpi) =\begin{cases} 
		-\sum_{\prA\subset \cA} \P_{\cA a}(\prA) Q(\pi_{\fb(\prA)}-\pi_a) \in[-1,0] &\text{ if }b=a,\\
		\sum_{\substack{\prA\subset\cA \\ \text{s.t. }b\in \fb(\prA)} } \P_{\cA a}(\prA) Q(\pi_{\fb(\prA)}-\pi_a) \in[0,1] &\text{ otherwise}
	\end{cases}$$
	for each action $b\in\cA$. Therefore,
	$$ \pdif{g_{a\fb}}{\bpi}(\x,\bpi) \Delta\bpi=\sum_{\prA\subset \cA} \P_{\cA a}(\prA) Q(\pi_{\fb(\prA)}-\pi_a)(\Delta\pi_{\fb(\prA)}-\Delta\pi_a) \quad\text{ for any }\Delta\bpi\in\R^\cA.$$
	
	As $g_{a\fb}$ is Lipschitz continuous on $\R^\cA$, the maximal value function $g_{a\ast}$ is also Lipschitz continuous and thus differentiable almost everywhere in $\R^\cA$. Furthermore, we have $\tpdif{g_{a\fb}}{\bpi}=\tpdif{g_{a\ast}}{\bpi}$ with $\fb=b_\ast$ at these differentiable points by applying a version of Danskin's envelop theorem (H\&S, Theorem A.4) to partial derivatives. By applying the above equation to $\fb=\fb_\ast$ and noticing $\pi_{\fb_\ast(\prA)}=\pi_\ast(\prA)$, we obtain
	$$ \dif{g_{a\ast}}{\bpi} (\x,\bpi) \Delta\bpi=\sum_{\prA\subset \cA} \P_{\cA a}(\prA) Q(\breve\pi_{a\ast}(\prA)) (\Delta\pi_\ast(\prA)-\Delta\pi_a) \quad\text{ for any }\Delta\bpi\in\R^\cA. $$
	Note that $g_{a\ast}$ is differentiable at $\bpi$ if and only if $\P_Q$ is continuous at $\breve\pi_{a\ast}(\prA)$ for all $a\in\cA$ and all $\prA\subset\cA$ with $\P_{\cA a}(\prA)>0$. In this case, any $\z_a \in\V_a(\x,\bpi)$ satisfies
	$$\z_a \cdot\Delta\bpi = \sum_{\prA\subset \cA} \P_{\cA a}(\prA) Q(\breve\pi_{a\ast}(\prA)) (\y_a(\prA)-\e_a)\cdot \Delta\bpi. 
	$$
	with some $\y_a(\prA)\in\Delta^\cA(b_\ast(\bpi,\prA))$. We have $\y_a(\prA)\cdot \Delta\bpi=\Delta\pi_\ast(\prA)$ and $\e_a\cdot \Delta\bpi=\Delta\pi_a.$ 
	Combining these, we obtain \eqref{eq:Dg_pi}.
\end{proof}

\subsubsection*{Proof of \Cref{clm:g_pi}, part a)}
First of all, we rigorously define $p^i_{\cA a}$ in property \ref{prop:PIndp_a}. Given payoff vector $\bpi$, make a partition of $\cA$ according to $\bpi$, say $\cA_1,\cA_2,\ldots$ such as 
$$ [a,a'\in \cA_i \ \Leftrightarrow \pi_a=\pi_{a'}] \qandq [a\in \cA_i \text{ and } a'\in \cA_{i'} \text{ with }i<i' \ \Leftrightarrow\ \pi_a>\pi_{a'}].$$
Let $\pi_i$ be the payoff obtained from actions in set $\cA_i$: i.e., $\pi_i$ is defined as $\pi_a$ with any $a\in \cA_i$. In other words, $\cA_i$ is the set of actions that yield the $i$-th greatest payoff $\pi_i$. Then, define $p_{\cA a}^i$ by
\begin{equation}
	p_{\cA a}^i \coloneqq \P_{\cA a}(\{\prA  \mid \prA\cap \cA_i\ne\emptyset \text{ and } \prA\cap \cA_j=\emptyset \text{ for any }j<i\}) \qquad\text{ for each } i.\label{eq:availIth}
\end{equation}
Given current action $a$, the maximal available payoff after the revision is $\pi_i$ with probability $p_{\cA a}^i$. For the set $\prA$ in the above definition of $p_{\cA a}^i$, we have $\pi_\ast(\prA)=\pi_i$. %

With this on hand, we can give a mathematical representation of property \ref{prop:PIndp_a}. Denote by $I(a)$ the index of the partition that $a$ belongs to: i.e., $a\in \cA_{I(a)}$. Then, property \ref{prop:PIndp_a} means
	\begin{equation}
	p_{\cA a}^i= p_{\cA b}^i \quad\text{ for any }i< \min\{I(a),I(b)\}. \label{eq:ASymmetry_Rank}
\end{equation}
Then we are now ready to prove part a) of \Cref{clm:g_pi}; for now, we postpone the proof of part b) since it involves a property of the aggregate net gain $G$.

\begin{proof}[Proof of \Cref{clm:g_pi}, part a)]
	\noindent\textbf{\Cref{assmp:AIndp_a_gen} $\Rightarrow$ property \ref{prop:PIndp_a}.} 
	If $\pi_a\le\pi_b$, then $I(a)\ge I(b)$. Fix $I< I(b)$ arbitrarily. Since $a,b\notin\bigcup_{i=1}^I \cA_i$ as $I< I(b)\le I(a)$, we obtain 
	$$ \sum_{i=1}^I p_{\cA a}^i = \P_{\cA a} \left(\{\prA \mid \prA \cap (\bigcup_{i=1}^I \cA_i) \ne \emptyset\} \right) = \P_{\cA b}\left(\{\prA_b \mid \prA_b \cap (\bigcup_{i=1}^I \cA_i) \ne \emptyset\} \right) = \sum_{i=1}^I p_{\cA b}^i $$
	for all $I< I(b)$ by applying \Cref{assmp:AIndp_a_gen} to $\cA_{ab}=\bigcup_{i=1}^I \cA_i$. Hence, we obtain \eqref{eq:ASymmetry_Rank}.
	
	\noindent\textbf{Property \ref{prop:PIndp_a} $\Rightarrow$ property \ref{prop:g_pi}.} Notice that 
	$$ g_{a\ast}(\x,\bpi)= \sum_{i=1}^{I(a)-1}p_{\cA a}^i \E_Q[\pi_i-\pi_a-q]_+.$$ 
	Assume $\pi_a\le\pi_b$ again. This implies $I(a)\ge I(b)$ and
		\begin{equation}
		\E_Q[\pi_i-\pi_a-q]_+\ge\E_Q[\pi_i-\pi_b-q]_+ \text{ for any }i.\label{eq:Prf_g1_CondlExpGain} 
	\end{equation}
	 Combining it with \eqref{eq:ASymmetry_Rank}, we obtain part i) of property \ref{prop:g_pi}:
	\begin{align*}
		g_{a\ast}(\x,\bpi)&= \sum_{i=1}^{I(b)-1}p_{\cA a}^i \E_Q[\pi_i-\pi_a-q]_+ +\sum_{i=I(b)}^{I(a)-1}p_{\cA a}^i \E_Q[\pi_i-\pi_a-q]_+\\
		\ge g_{b\ast}(\x,\bpi)&=\sum_{i=1}^{I(b)-1}p_{\cA b}^i \E_Q[\pi_i-\pi_b-q]_+.
	\end{align*}
	
	Further, if $\pi_a<\pi_b$, then $I(a)>I(b)$. Under the assumption of $\0\in\V_a(\x,\bpi)$, we have $g_{a\ast}(\x,\bpi)>0$ by property \ref{prop:gZero}-ii). Thus, there exists a rank $i<I(a)$ such that $p_{\cA a}^i>0$ and $\E_Q[\pi_i-\pi_a-q]_+ >0$. If $i\ge I(b)$, then it implies a strict inequality in the equation above. Consider a case of $i<I(b)$. Note that $\E_Q[\pi_i-\pi_a-q]_+ >0$ implies $Q(\pi_i-\pi_a)>0$ and thus $\E_Q[\pi_i-\pi_a-q]_+>\E_Q[\pi_i-\pi_b-q]_+$. With the fact $p_{\cA a}^i>0$, it also implies a strict inequality in the equation above. So, anyway, we obtain $g_{a\ast}(\x,\bpi)>g_{b\ast}(\x,\bpi)$, i.e., part ii) of property \ref{prop:g_pi}.
	
	\noindent\textbf{Property \ref{prop:g_pi}  $\Rightarrow$ \Cref{assmp:AIndp_a_gen}.} 
	We prove the contrapositive. For this, we assume the negation of \Cref{assmp:AIndp_a_gen}, especially the existence of $a,b\in\cA$ and set $\cA_{ab}^\sharp \subset \cA\setminus\{a,b\}$ such that 
	\begin{equation}
		\P_{\cA a} (\{ \prA \mid  \prA\cap \cA_{ab}^\sharp  \ne \emptyset \}) < \P_{\cA b} (\{ \prA_b \mid  \prA_b\cap \cA_{ab}^\sharp  \ne \emptyset \}).  \label{eq:neg_AInv_gen}
	\end{equation}
	Below we construct the payoff vector $\bpi^\sharp$ that negates property \ref{prop:PIndp_a}. Specifically, $\pi^\sharp_a<\pi^\sharp_b$ but $g_{a\ast}(\bpi^\sharp) < g_{b\ast}(\bpi^\sharp)$. Pick $\bar q>0$ such that $Q(\bar q)>0$.  Let $\pi^\sharp_i =\bar q$ for all $i\in \cA_{ab}^\sharp $; $\pi^\sharp_b =e \in (0,\bar q)$, $\pi^\sharp_a =0$ and $\pi^\sharp_j =0.5e$ for all other actions $j$.\footnote{Following the notation in the earlier part of this proof, we have $p_{\cA a}^1<p_{\cA b}^1$ by \eqref{eq:neg_AInv_gen}. Thus, the negation of property \ref{prop:PIndp_a} is obtained.} We have
	\begin{align*}
		g_{a\ast}(\bpi^\sharp)&= p_{\cA a}^1 \E_Q[\bar q-0-q]_+ + p_{\cA a}^2 \E_Q[e-0-q]_+ + p_{\cA a}^3 \E_Q[0.5e-0-q]_+,\\
		g_{b\ast}(\bpi^\sharp)&= p_{\cA b}^1 \E_Q[\bar q-e-q]_+.
	\end{align*}
	Note that $g_{a\ast}(\bpi^\sharp)$ and $g_{b\ast}(\bpi^\sharp)$ are continuous in $e$. As $e$ approaches to 0,  their limits satisfy
	\begin{align*}
		\lim_{e\to 0} g_{a\ast}(\bpi)= p_{\cA a}^1 \E_Q[\bar q-q]_+ \quad \text{ and }\quad
		\lim_{e\to 0} g_{b\ast}(\bpi)= p_{\cA b}^1 \E_Q[\bar q-q]_+.
	\end{align*}
	Note that $Q(\bar q)>0$ implies $\E_Q[\bar q-q]_+>0.$ Since  $p_{\cA a}^1<p_{\cA b}^1$, we have $\lim_{e\to 0} g_{a\ast}(\bpi^\sharp)<\lim_{e\to 0} g_{b\ast}(\bpi^\sharp)$. Therefore, for a sufficiently small $e$, we have $  g_{a\ast}(\bpi^\sharp)< g_{b\ast}(\bpi^\sharp)$ despite $\pi^\sharp_b =e>0=\pi^\sharp_a$. With such an $e$, we define  $\bpi^\sharp$ as above and then we obtain the negation of property \ref{prop:g_pi}, especially its part i).\footnote{Note that, if we additionally assume \Cref{assmp:AExAnteavailable_gen}, we have $p_{\cA a}^1>0$ and thus $g_{a\ast}(\bpi^\sharp)>0$; so, this proof also negates part ii) of property \ref{prop:g_pi}.} Therefore, the contrapositive is proven.
\end{proof}

\subsubsection*{Individual second-order ex-ante net gain $\h_\ast$}

For property \ref{prop:gzNeg},  we identify the change in an agent's (first-order) ex-ante net gain by a switch of the agent's action. So, we define the \textbf{individual ex-ante second-order net gain} for an action-$a$ player as
\begin{equation}
	h_{a\ast}(\x,\bpi)\coloneqq   \sum_{\prA\subset \cA} \P_{\cA a}(\prA;\x) Q(\breve\pi_{a\ast}(\prA)) (g_{\ast\ast}(\x,\bpi;\prA) -g_{a\ast}(\x,\bpi)). \label{eq:h_dfn}
\end{equation}

The term $\P_{\cA a}(\prA;\x) Q(\breve\pi_{a\ast}(\prA))$ is the (greatest) probability\footnote{\label{ftnt:2ndNetGain_IndetQ}If $\P_Q$ is not continuous, $Q(\breve\pi_{a\ast}(\prA))$ is only the \textit{greatest} possible probability of switch as the switching probability can be any value in the range $\cQ(\breve\pi_{a\ast}(\prA))$. Yet, we keep this definition to define $h_{a\ast}$ as a (single-valued) function of $\bpi$. So $h_{a\ast}$ is lower semicontinuous thanks to upper semicontinuity of cumulative distribution function.} that, upon the receipt of the first revision opportunity, the agent receives available action set $\prA$ and chooses to switch from the current action $a$. If $Q(\breve\pi_{a\ast}(\prA))>0,$ i.e., an agent chooses switching from $a$ with a positive probability, then it must be the case that $\pi_\ast(\prA)\ge\pi_a$. 
Property \ref{prop:g_pi} implies that this switch changes the agent's first-order net gain from $g_{a\ast}(\x,\bpi)$ to $g_{\ast\ast}(\x,\bpi;\prA).$ Hence $h_{a\ast}(\x,\bpi)$ is the expected change in the first-order gain upon the first revision opportunity, as stated in property \ref{prop:h_gz} below. Since the switch to a better reply decreases the first-order net gain, part i) of property \ref{prop:hZero} is implied by the above formula. Its part ii) needs a bit more careful examination of part ii) of property \ref{prop:g_pi}.

\begin{lem}\label{clm:h}
	Assume that the individual ex-ante gain function $\g_\ast:\R^\cA\to\R^\cA$ satisfies properties \ref{prop:gZero,prop:g_pi}.\footnote{No other assumptions or properties are needed. Once property \ref{prop:g_pi}. is satisfied, any assumption is not needed.}
	Then, the individual ex-ante second-order gain function $\h_\ast=(h_{a\ast})_{a\in\cA}:\R^\cA\to\R^\cA$ satisfies the following properties.
	\begin{description}
		\item[h\label{prop:hZero}] $h_{a\ast}$ is lower semicontinuous. i) $h_{a\ast}(\x,\bpi)\le 0$ for any $a\in\cA$ and $\bpi\in\R^\cA$. ii) $h_{a\ast}(\x,\bpi)=0$ if and only if $\0\in \V_a(\x,\bpi)$.
		\item[gh\label{prop:h_gz}] For any $a\in\cA$ and $\bpi\in\R^\cA$, if $g_a\ast$ is differentiable at $\bpi$, we have
		$$ h_{a\ast}(\x,\bpi)= \z_a\cdot \g_\ast(\x,\bpi) \qquad\text{ for any }\z_a\in\V_a(\x,\bpi).$$
		If it is not differentiable, we have $ h_{a\ast}(\x,\bpi)\ge \z_a\cdot \g_\ast(\x,\bpi)$ instead.
	\end{description}
\end{lem}

\begin{proof}[Proof of \Cref{clm:h}]
	\noindent\textbf{\ref{prop:hZero}.} First of all, since $Q$ is upper semicontinuous and $\breve\pi_{a\ast}(\prA):=\pi_\ast(\prA)-\pi_a$ is a continuous function of $\bpi$, the composite $ Q(\breve\pi_{a\ast}(\prA))$ is upper semicontinuous in $\bpi$. Since  $g_{\ast\ast}(\x,\bpi;\prA) - g_{a\ast}(\x,\bpi)$ is continuous and non-positive, this implies lower semicontinuity of $h_{a\ast}$.
	
	Part i) of property \ref{prop:hZero} is immediate from the definition of $g_{\ast\ast}$. For part ii) of property \ref{prop:hZero} , \eqref{eq:h_dfn} implies that $h_{a\ast}(\x,\bpi)=0$ is equivalent to $\P_{\cA a}(\prA)>0 \Rightarrow \left[ Q(\breve\pi_{a\ast}(\prA))=0\text{ or }g_{\ast\ast}(\x,\bpi;\prA) = g_{a\ast}(\x,\bpi) \right]$. On the other hand, $\0\in \V_a(\x,\bpi)$ is equivalent to $g_{a\ast}(\x,\bpi)=0$ by part ii) of prop \ref{prop:gZero} and thus to $\P_{\cA a}(\prA)>0 \Rightarrow \E_Q[\breve\pi_{a\ast}(\prA)-q]_+ =0.$ The latter is equivalent to $Q(\breve\pi_{a\ast}(\prA))=0$. 
	
	So, $\0\in \V_a(\x,\bpi)$ is equivalent to $h_{a\ast}(\x,\bpi)=0$. Now assume $h_{a\ast}(\x,\bpi)=0$. If $Q(\breve\pi_{a\ast}(\prA))=0$ whenever $\P_{\cA a}(\prA)>0$, then we have $\0\in \V_a(\x,\bpi)$. Consider the case in which there exists $\prA$ such that $\P_{\cA a}(\prA)>0$, $Q(\breve\pi_{a\ast}(\prA))>0$ and $g_{\ast\ast}(\x,\bpi;\prA) = g_{a\ast}(\x,\bpi)$. By part ii) of property \ref{prop:g_pi}, the last condition implies $\pi_\ast(\prA)=\pi_a$ or $\0\in \V_a(\x,\bpi)$. If the former implies $Q(\breve\pi_{a\ast}(\prA))=0$ and thus a contradiction with the assumption of $Q(\breve\pi_{a\ast}(\prA))>0$ in the current case. So, anyway, we obtain $\0\in \V_a(\x,\bpi)$. 

	\noindent\textbf{\ref{prop:h_gz}} First of all, notice that differentiablity of $g_a\ast$ at $\bpi$ implies continuity of $\P_Q$ at $\breve\pi_{a\ast}(\prA)$ for all $a\in\cA$ and all $\prA\subset\cA$ with $\P_{\cA a}(\prA)>0$. Any optimal actions in $b_\ast(\bpi,\prA)$ yield the greatest payoff $\pi_\ast(\prA)$ among all the available actions in $\prA$; property \ref{prop:g_pi} implies that they yield the smallest first-order gain $g_{\ast\ast}(\x,\bpi;\prA)$ among those available actions. Thus, any mixture of those optimal actions $\y_a(\prA)\in\Delta^\cA(b_\ast(\bpi,\prA))$ satisfies
	$$ \y_a(\prA)\cdot \g_\ast(\x,\bpi) = g_{\ast\ast}(\x,\bpi;\prA).$$
	By the fact $\e_a\cdot \g_\ast(\x,\bpi) =g_{a\ast}(\x,\bpi)$ and the assumption $\z_a\in\V_a(\x,\bpi)$, we find from \eqref{eq:h_dfn} that this implies the equation in property \ref{prop:h_gz}. 
\end{proof}

\subsubsection*{Aggregate net gain $G$: proofs of \Cref{clm:GZero,clm:DG} and part b) of \Cref{clm:g_pi}}

Now we derive properties of aggregate net gain function $G$ from those of individual ex-ante net gain function $\g_\ast$. To link aggregate and individual gains, we need the link between $\V$ and $\V_a$ as below; it leads to the proof of \Cref{clm:GZero}.

\begin{lem}\label{clm:Vdecomp}
	Fix $\x\in\Delta^\cA$ and $\bpi\in\R^\cA$ arbitrarily. We have $\0\in\V(\x,\bpi)$ if and only if $\left[ x_a>0 \ \Rightarrow\ \0\in\V_a(\x,\bpi) \right]$ for each $a\in\cA$.
\end{lem}

\begin{proof}[Proof of \Cref{clm:Vdecomp}]
The ``if'' part would be immediate by the definition of $\V$. For the ``only-if" part, assume that there are actions $a$ such that $x_a>0$ and $\0\notin\V_a(\bpi,\x)$, i.e., the sources of outflows of agents. Let $\cA_0$ be the set of actions that yield the smallest payoff in $\bpi$ among those actions. The statement $\0\notin\V_a(\bpi,\x)$ means that, for each of $a\in\cA_0$, there is an availability set $\prA$ such that $\P_{\cA a}(\prA)>0$ and $\breve\pi_{a\ast}(\prA)>0$. Thus, there is an action $b\in\prA$ that yields a strictly greater payoff than $a$ or any actions in $\cA_0$ and agents may switch from $a$ to this $b$; so, there is an outflow from $\cA_0$ (to such $b$ for each $a\in A$). Further, since the actions in $\cA_0$ yields smaller payoffs than any other sources of outflows. Since agents switch to better actions, these outflows never come to these actions $\cA_0$. Thus, the actions in $\cA_0$ only decrease the mass of players of these actions. So, $\V(\x,\bpi)$ cannot be a rest point.
\end{proof}

\begin{proof}[Proof of \Cref{clm:GZero} (property \ref{prop:GZero}) ]
	Part i) of property \ref{prop:GZero} is immediately obtained from part i) of property \ref{prop:gZero} and the fact $\x\in\R^\cA_+$. For part ii), by part i) of property \ref{prop:gZero}, $G(\x,\bpi)=0$ is equivalent to $\left[ x_a>0 \ \Rightarrow\ g_{a\ast}(\x,\bpi)=0 \right]$ since $G(\x,\bpi)\coloneqq \sum x_a g_{a\ast}(\x,\bpi)$. The equivalence between $\0\in\V_a(\x,\bpi)$ and $g_{a\ast}(\x,\bpi)=0$ has been established in part ii) in property \ref{prop:gZero}. Thus, combining these equivalences with \Cref{clm:Vdecomp}, we obtain part ii) of property \ref{prop:GZero}.
\end{proof}

Now we proceed to the proof of \Cref{clm:DG}. To prepare for property \ref{prop:gzNeg}, we define the \textbf{aggregate second-order gain function} $H:\Delta^\cA\times\R^\cA \to\R$ as $H(\x,\bpi) \coloneqq \x\cdot \h_\ast(\x,\bpi). $  The following lemma leads to part b) of \Cref{clm:g_pi} and then property \ref{prop:gzNeg} in \Cref{clm:DG}. 

\begin{lem}\label{clm:H}
	Consider a cost-benefit rationalizable dynamic $\V$. Assume properties \ref{prop:hZero} and \ref{prop:h_gz}. Then, the aggregate first-order and second-order gain functions $G,H:\Delta^\cA\times\R^\cA \to\R$ satisfy property \ref{prop:gzNeg}. 
\end{lem}

\begin{proof}[Proof of \Cref{clm:H}]
	\noindent\textbf{a) of property \ref{prop:gzNeg}} It can be readily obtained from property \ref{prop:hZero}, similarly to the proof of property \ref{prop:GZero} in \Cref{clm:GZero}.
	
	\noindent\textbf{b) of property \ref{prop:gzNeg}} Since $\Delta\x\in \V(\x,\bpi)=\sum_{a\in\cA} x_a \V_a(\x,\bpi)$  by \eqref{eq:TransitSet}, there exists $(\z_a)_{a\in\cA}$ such that $\Delta\x=\sum_{a\in\cA}x_a \z_a$ and $\z_a\in\V_a(\x,\bpi)$ for each $a$. By property \ref{prop:h_gz}, we have
	$$ \Delta\x\cdot\g_\ast(\x,\bpi) =\sum_{a\in\cA} x_a \z_a\cdot\g_\ast(\x,\bpi)=\sum_{a\in\cA} x_a h_{a\ast}(\x,\bpi)=H(\x,\bpi).\vspace{-2\baselineskip} $$
\end{proof}

\begin{proof}[Proof of \Cref{clm:g_pi}, part b)]
	Property \ref{prop:g_pi} implies properties \ref{prop:hZero} and \ref{prop:h_gz} by \Cref{clm:h}, which are succeeded to properties \ref{prop:gzNeg} according to \Cref{clm:H}. 
\end{proof}
	
\begin{proof}[Proof of \Cref{clm:DG}]
	\noindent\textbf{\ref{prop:GDpi}.} As $G(\x,\bpi)=\x\cdot\g_\ast(\x,\bpi)$, property \ref{prop:gDpi} in \Cref{clm:g0} implies
	$$ \x\cdot \pdif{\g_\ast}{\bpi}(\x,\bpi) \Delta\bpi =\sum_{a\in\cA} x_a \pdif{g_{a\ast}}{\bpi}(\x,\bpi)\Delta\bpi =\sum_{a\in\cA} x_a\z_a \cdot\Delta\bpi$$
	for any $\z_a\in\V_a(\x,\bpi)$. Since $\Delta\x\in\V(\x,\bpi)=\sum_{a\in\cA} x_a\V_a(\x,\bpi)$ by \eqref{eq:TransitSet}, we obtain property \ref{prop:GDpi}. 
	
	\noindent\textbf{\ref{prop:gzNeg}.} \Cref{clm:g_pi} suggests that \Cref{assmp:AIndp_a_gen} is equivalent to property \ref{prop:g_pi} by part a) and it further implies properties \ref{prop:gzNeg} by part b).
	
	\noindent\textbf{\ref{prop:gDx}.}  This is immediate from \Cref{assmp:AIndp_x_gen}, since $\g_\ast$ could vary with $\x$ only through dependency of $\P_{\cA a}$ on $\x$ according to the definition of $g_{a\ast}$ as in \eqref{eq:g_dfn}. 
\end{proof}

\subsection{Proofs of other theorems}
\subsubsection*{Proof of \Cref{clm:StatDynStbl}}
\begin{proof}
		We use function $G$ as a Lyapunov function $W$ in \Cref{thm:Lyapunov_DI} and function $H$ as a decay rate function $\tilde W$. Thanks to continuous differentiability (and thus Lipschitz continuity) of $\F$, Lipschitz continuity of $G$ and lower semicontinuity of $H$ are succeeded to those of $G^\F$ and $H^\F$, respectively; condition (c) in the theorem is met. Properties \ref{prop:GZero} and \ref{prop:gzNeg} (part a) meet conditions (a) and (b). 
			
	By calculating $\dot G^F$ with the chain rule and then applying properties \ref{prop:gDpi}, \ref{prop:gzNeg} (part b)  and \ref{prop:gDx}, we have 
	\begin{align}
		\dif{}{t}G^\F(\x_t)&\equiv \pdif{G}{\x}(\x_t,\F(\x_t))\dot\x_t +\pdif{G}{\bpi}(\x_t,\F(\x_t))\dif{\F}{\x}(\x_t)\dot\x_t \notag\\
		&=\dot\x_t \cdot\dif{\F}{\x}(\x_t)\dot\x_t +H^\F(\x_t) \label{eq:dGdt_H}
	\end{align}
	for any $\dot\x_t\in\V^\F(\x_t)$ and $\x_t\in\Delta^\cA$. If $D\F$ is negative semidefinite at $\x_t$, \eqref{eq:dGdt_H} implies
	$$ \dif{}{t}G^\F(\x_t) \le H^\F(\x_t)\le 0\qquad\text{ for any $\dot\x_t\in\V^\F(\x)$}.$$
	Thus we obtain \eqref{eq:LCondtn_DI} without making restriction on the initial state $\x_0$; so it holds globally and thus the basin of attraction is $\Delta^\cA$. This proves part i). By the same token, we can verify the local stability as in part ii-a) by setting $X'$ to a neighborhood of $\x^\ast$ in which $\x^\ast$ is the only rest point and $D\F$ is negative semidefinite.
	
	For part ii-b), consider a Carath\'{e}odory solution path $\{\x_t\}_{t\in\R_+}$ from $\x_0=\x^\ast$. Then, we have $\dot G^\F(\x_0)=H^\F(\x^\ast)+\dot\x_0\cdot D\F(\x^\ast)\dot\x_0$. The former term $H^\F(\x^\ast)$ is zero by part a) of property \ref{prop:gzNeg}, since $\x^\ast$ is a rest point. The assumption for part ii-b) allows us to take $\dot\x_0\in\V^\F(\x^\ast)$ that makes the latter term $\dot\x_0\cdot D\F(\x^\ast)\dot\x_0$ being positive. With $G(\x_0)\ge 0$ (indeed zero by part ii) of property \ref{prop:GZero}), this implies $G^\F(\x_t)>0$ for any $t$ in some interval of time $(0,T]$. Part ii) of property \ref{prop:GZero} implies that $\x_t$ in this interval of time cannot be at a rest point. Thus, it must be escaping from $\x^\ast$; this implies that $\x^\ast$ is not Lyapunov stable.   
\end{proof}

\subsubsection*{Proof of \Cref{clm:ReplRefuelG}}
First we prove an algebraic lemma as below.\footnote{In the lemma, $A$ can be any natural number and only means the dimension of vectors.}
\begin{lem}\label{clm:TripleProd}
	Consider three $A$-dimensional vectors $\x\coloneqq (x_1,\ldots,x_A), \y\coloneqq(y_1,\ldots,y_A), \w\coloneqq(w_1,\ldots,w_A)$ such that (assumption xy) $\sum_{a=1}^\cA x_a y_a = 0$, (x) $\x\in\Delta^\cA$, (y)  $y_1\le y_2 \le \cdots \le y_A$ and (w) $w_1\le w_2 \le \cdots \le w_A$. Then, we have the following.
	\begin{enumerate}[i)]
		\item $\sum_{a=1}^\cA x_a y_a w_a \ge 0$. 
		\item Furthermore, if (xyw) there exists indexes $b,c$ such that $x_b,x_c>0$ and $y_b<0<y_c$ and $w_b<w_c$, then we have $\sum_{a=1}^\cA x_a y_a w_a > 0$.
	\end{enumerate} 
\end{lem}
\begin{proof}[Proof of \Cref{clm:TripleProd}]
	\noindent\textbf{i).} Assumptions (xy), (x), and (y) imply the existence of indexes $I,J$ such that $1\le I<J\le A$ and $ y_a \le 0 $ for any $a\le I$, $y_a= 0$ for any $a\in (I,J),$ and $y_a\ge 0$ for any $a\ge J. $ Pick any $\bar w\in [w_I,w_J]$, noticing (w). Then, we have 
	\begin{equation}
			\sum_{a=1}^I x_a y_a (\bar w-w_a) \le \sum_{a=I+1}^{J-1} x_a y_a (w_a-\bar w)=0 \le \sum_{a=J}^A x_a y_a (w_a-\bar w). \label{eq:TripleProd_prf}
	\end{equation}
	since $x_a\ge 0$, $y_a\le 0$ and $w_a\le \bar w$ for any $a\le I$; $y_a=0$ for any $a\in [I+1,J-1]$ (if any); and, $x_a\ge 0$, $y_a\ge 0$ and $w_a\ge \bar w$ for any $a\ge J$. Thus, we have $  \sum_{a=1}^\cA x_a y_a (w_a-\bar w) \ge 0.$ By (xy), this implies
	$$ \sum_{a=1}^\cA x_a y_a w_a \ge  \bar w \sum_{a=1}^\cA x_a y_a = 0.$$
	\noindent\textbf{ii).} Assumption (xyw) implies that $w_b<\bar w$ or $w_c>\bar w$ holds for any choice of $\bar w\in [w_I,w_J]$, as well as $b\le I<J\le c$. Hence, $x_b y_b (\bar w-w_b)$ or $x_c y_c (w_c-\bar w)$ is strictly positive. These imply that the above equation holds with a strict inequality. 
\end{proof}

\begin{proof}[Proof of \Cref{clm:ReplRefuelG}]
	\noindent\textbf{i)} Given $\F(\x)$, we reorder actions by the reverse of payoff ordering so $F_i(\x)\le F_j(\x)$ for any $i< j$. Then, the availability effect in $\dot G^\F$ under the replicator dynamic can be written as 
	\begin{align*}
		G_\text{av} &\coloneqq \x \cdot \pdif{\g_\ast}{\x}(\x,\F(\x)) \dot\x 
			= \sum_{i\in\cA} x_i \sum_{j\in\cA} \dot x_j \E_Q[F_j(\x)-F_i(\x)-q]_+ \\
			&= \sum_{j\in \cA} x_j F_j(\x) \sum_{i\le j}  x_i \E_Q[F_j(\x)-F_i(\x)-q]_+.
	\end{align*}
	The second equation comes from $g_{i\ast}(\x,\bpi)=\sum_j x_j  \E_Q[\pi_j-\pi_i-q]_+$ in the replicator dynamic ($\P_\cA(\{ j \})=x_j$) and thus  $\tpdif{g_{i\ast}}{x_j}=\E_Q[\pi_j-\pi_i-q]_+$. The third equation comes from $\dot x_j=x_j (F_j(\x)-\x\cdot\F(\x))$ in the replicator dynamic and $\x\cdot\F(\x)\equiv 0$ in a zero-sum game; note that, if $i>j$, we have $F_i(\x)\ge F_j(\x)$ and thus $F_j(\x)-F_i(\x)-q\le 0$ for any realization of $q\in\R_+$. 
	
	We apply \Cref{clm:TripleProd} to this $\x$ with $\y=\F(\x)$ and $\w$ such that $w_j=\sum_{i\le j}  x_i \E_Q[F_j(\x)-F_i(\x)-q]_+$. (Note that $A$ in the theorem is $\sharp\cA$ and $w_A=0.$) Then, they satisfy the assumptions for part i) of \Cref{clm:TripleProd}; note that we have $\sum x_a y_a=\sum x_a F_a(\x)=0$ since $\F$ is a zero-sum game; and, $w_a\ge w_b$ whenever $a\ge b$ since $x_i\E_Q[F_a(\x)-F_i(\x)-q]_+\ge x_i\E_Q[F_b(\x)-F_i(\x)-q]_+\ge 0$ for any $i\le b(\le a)$ by $F_a(\x)\ge F_b(\x)$ and $x_i \E_Q[F_a(\x)-F_i(\x)-q]_+\ge 0$ for any $i\in [b,a)$. Thus, we have $G_\text{av}\ge 0$. 
	
	\noindent\textbf{ii)}  If $\x$ is not a rest point under the replicator dynamic in a zero-sum game, then there must be an action $b$ such that $x_b>0$ and $y_b=F_b(\x)\ne \x\cdot\F(\x)=0$.\footnote{If $\x$ is a rest point under the replicator dynamic, then each $a$ satisfies $x_a=0$ or $F_a(\x)=\x\cdot\F(\x)$.} Besides, because $\x\cdot\F(\x)=0$, if $y_b<0$, there is another action $c$ such that $x_c>0$ and $y_c=F_c(\x)>0$. Thus, without loss of generality, let $y_b<0<y_c$; note that $b<c$. Since $Q(\bar q)=[\bar q]_+$ and thus $\E_Q[\bar q-q]_+ = 0.5[\bar q]_+^2$ in the replicator dynamics, we have
	\begin{align*}
		w_c - w_b & =  \sum_{i\le c}  x_i \E_Q[F_c(\x)-F_i(\x)-q]_+ - \sum_{i\le b}  x_i \E_Q[F_b(\x)-F_i(\x)-q]_+  \\
			& \ge \sum_{i\le b}  x_i \left(\E_Q[F_c(\x)-F_i(\x)-q]_+ - \E_Q[F_b(\x)-F_i(\x)-q]_+\right) \\
			&\ge x_b \left(\E_Q[F_c(\x)-F_b(\x)-q]_+ - \E_Q[F_b(\x)-F_b(\x)-q]_+\right) \\
			&=  0.5 x_b (F_c(\x)-F_b(\x))^2 >0.
	\end{align*}
	Hence, assumption (xyw) is satisfied. Thus, part ii) of  of \Cref{clm:TripleProd} is applicable and we obtain $G_\text{av}>0$.
\end{proof}

\section{Appendix to \Cref{sec:MainResults}\label{apdx:MainResults}}
\subsection{About \Cref{assmp:AIndp_a_gen}\label{apdx:assumpAInv}}
\subsubsection*{Why not a simpler definition?}
The role of \Cref{assmp:AIndp_a_gen} in \Cref{clm:StatDynStbl} is to obtain property \ref{prop:PIndp_a}, which is equivalent to property \ref{prop:g_pi} and thus leads to property \ref{prop:gzNeg}. \Cref{assmp:AIndp_a_gen} implies that, given the preference ordering of actions based on the current payoff vector, for any $i\in\N$, the probability that (at least one of) the $i$-th best actions is available does not vary whether the current action is $a$ or $b$, as long as the $i$-th action is better than both $a$ and $b$. This is what property \ref{prop:PIndp_a} means. So, if $a$ is worse than $b$, an agent who is currently taking $a$ is more likely to switch than one currently taking $b$ because the former action $a$ has more actions (including $b$) that are better than $a$ itself and thus the agent should switch to. Besides, if the former agent and the latter switch to the same action, the former improves the payoff more than the latter, as in \eqref{eq:Prf_g1_CondlExpGain}. Combining these two facts, the ex-ante gain by a switch from $a$ should be greater than that from $b$. Therefore, the ordering of actions by the ex-ante gains is obtained by simply reversing the ordering by the current payoffs, as stated in property \ref{prop:g_pi}.\footnote{The switching cost does not matter to this and thus the monotone relationship as in property \ref{prop:g_pi} is obtained as well for the gross gain.}

\Cref{assmp:AIndp_a_gen} may look a bit complicated. First, it checks the symmetry for each \textit{set }$\cA_{ab}\subset \cA\setminus\{a,b\}$, not for each action $c\in\cA\setminus\{a,b\}$. Second, it checks the probability of the event that \textit{at least one action }in the set is available, not simply that the whole set is available. These assure invariance of correlations in availability between any actions other than $a$ and $b$; particularly, it allows us to distinguish the availability of the second best alone from that of the first and second best together, which results in different amounts of payoff gains. 

One might think of somewhat simpler conditions that look similar to \Cref{assmp:AIndp_a_gen}:
\begin{enumerate}[i)]
	\item $\P_{\cA a}(\{\prA \mid c\in \prA\})=\P_{\cA b}(\{\prA \mid c\in \prA\})$ for any $a,b\in\cA$ and $c\in \cA\setminus\{a,b\}$. That is, the probability that \textit{arbitrary action} $c$ is available should not depend on the status quo, as long as strategy $c$ is not the status quo.
	\item $\P_{\cA a}(\cA_{ab})=\P_{\cA b}(\cA_{ab})$ for any $a,b\in\cA$ and $\cA_{ab}\subset \cA\setminus\{a,b\}$. That is, for any subset that excludes both $a$ and $b$, the probability that the available strategy set \textit{coincides exactly} with this subset does not depend on which of the excluded two actions is the status quo.
\end{enumerate}
However, neither of these two simpler conditions cannot assure property \ref{prop:PIndp_a}. Each of the following cases satisfies one of the above two conditions but not \Cref{assmp:AIndp_a_gen}; as a result, property \ref{prop:PIndp_a} does not hold. Here we assume $\cA=\{a,b,c,d,e\}$ and payoff vector $\bpi$ such that $\pi_e>\pi_d>\pi_c>\pi_b>\pi_a$. 

\paragraph{i)} Assume that, when $a$ is the status quo, then the available strategy set is $\{b\}$ with probability 0.4, $\{c,d\}$ with 0.3 and $\{e\}$ with 0.3. When $b$ is the status quo, then the available strategy set is $\{a\}$ with probability 0.4, $\{c\}$ with 0.3 and $\{d,e\}$ with 0.3.
Condition i) is satisfied, since each of the three actions $c,d$ and $e$ is available with probability 0.3, whether the status quo is $a$ or $b$. However, \Cref{assmp:AIndp_a_gen} is not satisfied since at least $d$ \textit{or }$e$ is available with probability $\P_{\cA a}(\{\prA \mid \prA\cap\{d,e\}\neq\emptyset\})=\P_{\cA a}(\{c,d\})+ \P_{\cA a}(\{e\})= 0.6$ when $a$ is the status quo, while this probability is $\P_{\cA b}(\{\prA \mid \prA\cap\{d,e\}\neq\emptyset\})=\P_{\cA b}(\{d,e\})=0.3$ when $b$ is the status quo. 

Note that, if the second best $d$ is available, the first best $e$ is not available when $a$ is the status quo; when $b$ is the status quo, $e$ is also available. As a result, the probability that the second best $d$ is the best available strategy is $p^2_{\cA a}=0.3$ when $a$ is the status quo, while it is $p^2_{\cA b}=0$ when $b$ is the status quo. Thus, property \ref{prop:PIndp_a}  does not hold. This is because condition i) allows correlation of availabilities of distinct two strategies to vary with the status quo. 

\paragraph{ii)} Assume that, when $a$ is the status quo, then the available strategy set is $\{b\}$ with probability 0.8, $\{b,c,d\}$ with 0.1, and $\{b,e\}$ with 0.1. When $b$ is the status quo, then the available strategy set is $\{a\}$ with probability 0.8, $\{a,c\}$ with 0.1, and $\{a,d,e\}$ with 0.1.
Condition ii) is trivially satisfied, since no subset that excludes both $a$ and $b$ is available with any positive probability when either $a$ or $b$ is the status quo. (Note that condition i) is also satisfied.) However, \Cref{assmp:AIndp_a_gen} is not satisfied, since at least $d$ \textit{or }$e$ is available with probability $0.2$ when $a$ is the status quo; this probability is $0.1$ when $b$ is the status quo. 

Similarly to the example for condition i), $\P_{\cA a}$ and $\P_{\cA b}$ have different correlations between availabilities of the first best action $e$ and the second best action $d$; this results in $p^2_{\cA a}=0.1\ne 0=p^2_{\cA b}$. Thus, property \ref{prop:PIndp_a} does not hold.

\subsubsection*{Non-trivial interesting examples}
\Cref{assmp:AIndp_a_gen} allows $\P_{\cA \cdot}$ to be different over different actions---some actions can be more likely to be available than others---as long as the availability does not depend on the current action. This is utilized in the following two examples; it is easy to confirm that they satisfy full rationalizability.
Besides, the generalized assumption illustrates the difference from payoff monotonicity; see \Cref{apdx:OrderComp}.

\begin{exmp}[Non-uniform availability]
	Assume that the action set $\cA$ is partitioned to $\cA_1,\ldots,\cA_I$ (with $I\in\N$).\footnote{As a partition of $\cA$, they satisfy $\cA=\cup_{i=1}^I \cA_i$ and $\cA_i\cap \cA_j=\emptyset$ if $i\ne j$. This partition is different from the one in the proof of \Cref{clm:g_pi}, part a) in \Cref{apdx:Proofs}.} An action in $\cA_i$ is available with probability $p_i$. Assume that availability of an action is independent of other actions' availability (including those of actions in the same set $\cA_i$) and that there are at least two actions that are always available; say, $\cA_1$ with $p_1=1$ and $\#\cA_1\ge 2$.\footnote{This is to guarantee $\P_{\cA a}(\{\prA \mid \prA \ne\emptyset\})=1$ for any $a$; that is, the available action set $\cA_a'$ is non-empty with probability 1 whatever the status-quo action $a$ is.} Then, we have
	$$ \P_{Aa}(\prA) =\prod_{i=1}^I p_i^{\#(\cA_i\cap \prA)} \qquad\text{ for any }\prA \subset \cA.\vspace{-2\baselineskip}$$
\end{exmp}

\begin{exmp}[Partitioned available action sets]
	We can have actions available only in a group. Again, consider a partition of $\cA$ to $\cA_1,\ldots,\cA_I$ and assume $\#\cA_i\ge 2$ for all $i$.\footnote{Without this assumption, we may have $\cA_i\setminus\{a\}=\emptyset$. Then, if the status-quo action is $a$, a revising agent may have no available action (except the status quo) with probability $p_i$.} Either one of those action subsets is randomly picked but not more than one; when $\cA_i$ is picked, all the actions in $\cA_i$ are available for the agent. Noticing that status-quo action $a$ should be excluded from the available action set (incurring a switching cost), we can formalize this procedure as
	$$ \P_{Aa}(\cA_i\setminus\{a\})=p_i \qquad \text{ with }\sum_{i=1}^I p_i=1.\vspace{-2\baselineskip}$$
\end{exmp}

\section{Appendix to \Cref{sec:Ext}\label{apdx:Ext}}
\subsection{Appendix to \Cref{sec:Modified}}\label{apdx:Modified}
\subsubsection*{Modified framework for excess payoff dynamics}
To interpret excess payoff dynamics as cost-benefit rationalizable, we need just one modification of our framework: the status-quo strategy should be $\x$ as a mixed strategy, i.e., the average strategy (action distribution) of actions taken by agents in the population. We imagine a birth-death process over generations \citep{AlosferrerNeustadt08IGTR_BRDinBirthDeath}, rather than a revision process of immortal agents as in the basic framework.\footnote{
	Specifically, an agent's death follows a Poisson process with arrival rate 1. At the same time as one's death, a new agent is born. Unless a newborn agent chooses to pay $q$, the agent is assigned to a randomly chosen action, drawn from the population's action distribution $\x$.} 
Everything else is the same as the basic framework. If an agent chooses to switch an action, then the agent chooses a pure strategy of a particular available action and pays switching cost $q$. Available action set $\prA\subset \cA$ is randomly drawn from probability distribution $\P_\cA$ over the power set of $\cA$; switching cost $q$ follows probability distribution $\P_Q$ over $\R_+$. A newborn agent compares the net payoffs from pure strategies of available actions with the expected payoff from the default mixed strategy $\x$; then the newborn agent decides whether to adopt the default strategy $\x$ or to pay cost $q$ and take a pure strategy of an available action.

The gross gain from a switch to action $a$ is indeed the excess payoff of the action, i.e., $\hat\pi_a\coloneqq \pi_a-\bpi\cdot\x$; let $\hat\pi_\ast(\prA)\coloneqq \pi_\ast(\prA)-\bpi\cdot\x$. Thus, after drawing available action set $\prA$, an agent's conditional switching rate must belong to $\cQ(\hat\pi_\ast(\prA))$. Thus, the transition of the social state follows
$$ \dot\x\in \V(\x,\bpi) =\sum_{\prA\subset \cA} \P_\cA(\prA;\x) \cQ(\hat\pi_\ast(\prA)) \left( \Delta^\cA(b_\ast(\prA,\bpi)-\x \right).$$

We call an evolutionary dynamic \textbf{a cost-benefit rationalizable birth-death dynamic} if it is constructed in this way.  As this construction of $\P_\cA$ meets \Cref{assmp:AIndp_a}, full rationalizability only requires \Cref{assmp:AIndp_x}.\footnote{However, we could have $\P_{\cA a}$ dependent on a revising agent's current action $a$. Possibly, if we still interpret it as a birth-death process, we can say that availability of actions $\P_{\cA a}$ for a newborn child still depends on the child's parent but the parent does not have an influence on the child's choice of an action; rather, the social average does.} 

This construction implies that a revising agent takes only actions that yield greater payoffs than the population average. This better reply property may not hold for integrable excess payoff dynamics.\footnote{This is indeed a disturbance to understand integrability economically. \cite{Sandholm05JET_ExcessPayDyn} or \cite{HofSand09JET_StableGames} does not give a concrete example of such a dynamic.} 

On the other hand, consider a cost-benefit rationalizable birth-death dynamic with $\P_\cA(\{a,b\})=2/\sharp A(\sharp A-1)$ for any pair of two distinct actions $a,b$ (satisfying \Cref{assmp:AIndp_x}, so we drop $\x$ from arguments of  $\P_\cA$): a revising agent faces a randomly picked pair of two actions as an available action set and then switches to a better action if it yields a greater payoff than the population average. Although the sampling process is different from the one that straightforwardly induces an excess payoff dynamic as in \Cref{fig:majorDyn}, this dynamic can be regarded as a version of excess payoff dynamics since the rate of switch to each action is determined from relative payoff vector $\hat\bpi$.\footnote{Note that this dynamic is not separable as the switching rate to one action depends on its payoff rank compared with other actions.} This dynamic exhibits discontinuity in the switching rate when two actions yielded an equal payoff before but either one becomes better than the other. Indeed, when there is a tie in payoffs, this dynamic allows for multiple transition vectors. By construction, this dynamic is fully rationalizable. Yet, since the integrability condition by H\&S (eq. (1)) requires continuity of the switching rate function and also their class of excess payoff dynamics is defined by a differential equation, this dynamic is not covered by a class of integrable dynamics.

\subsubsection*{Equilibrium stability}
Now we can define an individual's first-order gain $g_\ast(\x,\bpi)$ as 
$$ g_\ast(\x,\bpi)\coloneqq   \sum_{\prA\subset \cA} \P_\cA(\prA;\x) \E_Q[\hat\pi_\ast(\prA)-q]_+. $$
With the total mass of agents equal to one, the aggregate first-order gain $G(\x,\bpi)$ is now just equal to the individual first-order gain: $G(\x,\bpi)=g_\ast(\x,\bpi)\cdot 1$.
Similar to \Cref{clm:GZero,clm:DG}, we can prove properties \ref{prop:GZero} and \ref{prop:GDpi}. For availability effect, we define function $H$ by 
$$ H(\x,\bpi) \coloneqq -\left\{\sum_{\prA\subset \cA} \P_\cA(\prA;\x) Q(\hat\pi_\ast(\prA))\right\}\left\{\sum_{\prA\subset \cA} \P_\cA(\prA;\x) Q(\hat\pi_\ast(\prA)) \hat\pi_\ast(\prA)\right\}.$$
Under \Cref{assmp:AIndp_x}, the right hand side is indeed $\Delta\x\cdot \tpdif{\g_\ast}{\x}(\x,\bpi)$ and $H$ satisfies $H(\x,\bpi)$ is non-positive at any $(\x,\bpi)$ and equals to zero if and only if $\0\in\V(\x,\bpi)$. Furthermore, the switching effect is zero since $g_\ast$ does not depend on an agent's current action. In sum, we have 
$$  \dot G^\F(\x)= \x\cdot D\F(\x)\dot \x + H^\F(\x) +0.$$
Under stable stability of $\F$, this reduces to $\dot G^\F(x)\le H^\F(\x)$; so, $G^\F$ serves as a Lyapunov function in \Cref{thm:Lyapunov_DI} with $H^\F$ a decaying rate function.

\subsection{Appendix to \Cref{sec:hetero}}\label{apdx:hetero}
\subsubsection*{Set-up of a multi-population game.}
The society consists of $P$ populations; let $m^p>0$ be the mass of population $p\in\cP\coloneqq\{1,\ldots,P\}$ and $\sum_{p\in\cP} m^p=1$. Let $\cA^p=\{1,\ldots,A^p\}$ be the action set for an agent in population $p$. Denote $A^P:=\sum_{p}A^p$. We denote by $x^p_a$ the mass of action-$a$ agents in population $p$; population $p$'s state $\x^p\coloneqq (x^p_a)_{a\in\cA^p}$ belongs to $\cX^p\coloneqq m^p\Delta^{\cA^p}$. The state of the entire society is $\x^\cP\coloneqq (\x^p)_{p\in\cP}$. Denote by $\cX^\cP\coloneqq \bigtimes_{p\in\cP} \cX^p$ the set of all the feasible social states. Let $\F^p:\cX^\cP\to\R^{A^p}$ be the payoff function for an agent in population $p$. $\F^p(\x^\cP)$ is the payoff vector for agents in population $p\in\cP$ when the social state is $\x^\cP$. Let $\bpi^\cP$, given a collection of payoff vectors $\bpi^p\in\R^{A^p}$ over all the populations $p\in\cP$. The definition of static stability (\Cref{dfn:StatStbl}) is extended to a multi-population setting by requiring $\F,\x,\Delta\x$ in condition \eqref{eq:DF_StatStblCondn} to be replaced with $\F^\cP,\x^\cP,\Delta\x^\cP$. 

\begin{exmp}[Saddle games]
	It has been long to call game theorists' attention to a multi-player game with a ($C^2$-class) saddle function $\phi:\cX^\cP\to\R$ that is concave in $\x^C \coloneqq (\x^p)_{p\in\cP^C}$ and convex in $\x^V \coloneqq (\x^p)_{p\in\cP^V}$.\footnote{
		The study of evolutionary dynamics in this class of games dates back to \cite{Kose56EMA_SolSaddle_DiffEqm}.
		\cite{HofSorin06DCDS_BRD_ContZerosum} verify Nash stability in this class of games under the standard BRD, while they regard this game as a zero-sum game. A class of saddle games is broader than that of zero-sum games; see \citet[Example 1]{NoraUno14JET_SaddleFn_RobustEqia}.} Here, the population set is partitioned to $\cP^C$ and $\cP^V$.  The payoff function for population $p$ is given by $\F^c\equiv \tpdif{\phi}{\x^c}$ for each $c\in\cP^C$ and $\F^v\equiv -\tpdif{\phi}{\x^v}$ for each $v\in\cP^V.$ Then, $\F^\cP$ is a stable game and a Nash equilibrium $\x^{\cP\ast}$ is a saddle point of $\phi$. 
\end{exmp}

\begin{exmp}[Anonymous games]
	In an \textbf{anonymous game}, all populations $p\in\cP$ share the same action set $\cA$ and each population's payoff vector depends on $\x^\cP$ only through the aggregate action distribution over the society $\bar\x\coloneqq  \sum_{q\in\cP} \x^q \in \Delta^\cA$. Assume that payoff heterogeneity is additively separable: each $p$'s payoff function $\F^p$ is written as 	$\F^p(\x^\cP)=\F^0(\sum_{q\in\cP} \x^q)+\btheta^p$, the sum of common payoff function $\F^0:\Delta^\cA\to\R^\cA$ and constant payoff perturbation $\btheta^p\in\R^\cA$.
	If $\F^0$ is a stable game in the homogeneous setting, then $\F^\cP$ is also a stable game.
\end{exmp}

\subsubsection*{Evolutionary dynamics and gains}
Let $\V^p:\cX^p\times\R^{A^p}\to \R^{A^p}_0$ be the evolutionary dynamic for population $p$'s state $\x^p$. We allow each population $p$ to follow a different revision protocol but we assume that $\V^p$ depends only on $\x^p$ and $\bpi^p$.\footnote{If a dynamic satisfies \Cref{assmp:AIndp_x_gen}, it only requires uncoupledness---independence of the dynamic from $\bpi^{-p}$. Otherwise, it requires dependency of $\P_{\cA a}$ on the social state to be limited to $\x^p$ of the population that a revising agent belongs to.} Then, the social dynamic of $\x^\cP$ is constructed simply as 
$$ \dot\x^\cP=(\dot\x^p)_{p\in\cP} \in \V^\cP(\x^\cP,\bpi^\cP)\coloneqq \bigtimes_{p\in\cP} \V^p (\x^p,\bpi^p).$$
Suppose that each population's evolutionary dynamic $\V^p$ has functions $G^p,H^p:\cX^p\times\R^{A^p}\to\R$ that satisfies properties in \Cref{clm:GZero,clm:DG}. Then, we define functions $G^\cP,H^\cP:\cX^\cP\times \R^{A^\cP}\to\R$ for the entire society $\cP$ by $ G^\cP(\x^\cP,\bpi^\cP) \coloneqq \sum_{p\in\cP} G^p(\x^p,\bpi^p)$ and $H^\cP(\x^\cP,\bpi^\cP)\coloneqq \sum_{p\in\cP} H^p(\x^p,\bpi^p).$
As these are simply sums, it is easy to confirm that the properties of $G^p$ and $H^p$ are succeeded to those of $G^\cP$ and $H^\cP$.

Choices in different populations affect evolution in other populations through playing a game $\F^\cP$ on the same field.\footnote{Note that, because of this, an uncoupled dynamic cannot guarantee convergence to Nash equilibria in general, even if it satisfies \Cref{assmp:AIndp_x_gen} \citep{HartMascolell03AER_UncoupledDyn}.} However, separation of an evolutionary dynamic $\V(\x,\bpi)$ from a game $\bpi=\F(\x)$ allows us to put aside this indirect interdependence between  populations.

\section{Appendix to \Cref{sec:dis}}\label{apdx:dis}
\subsection{Appendix to \Cref{sec:GrossGain}}
The following example of a strictly stable game illustrates a case where the aggregate \textit{gross} gain may not decrease over time even under a cost-benefit rationalizable dynamic.%

\begin{exmp}[Rock-Paper-Scissors\label{exmp:RPS}]
	\begin{figure}[t!]
		\subfloat[Trajectory under Smith\label{fig:RPS_Smith_Str}]{\includegraphics[width=0.44\textwidth]{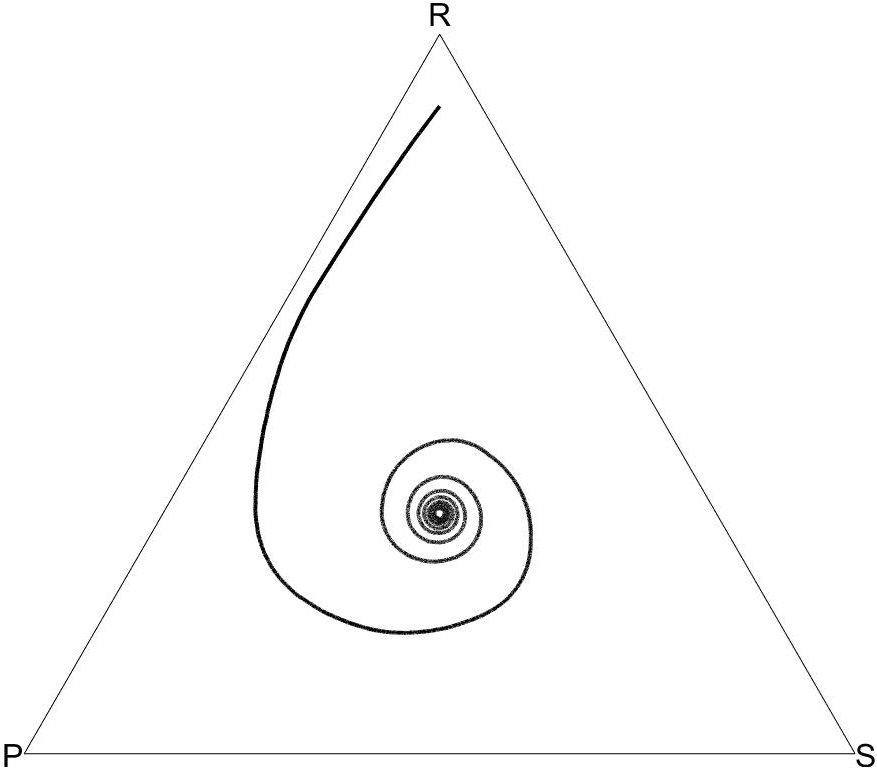}}\hfill
		\subfloat[Net and aggregate gains under Smith\label{fig:RPS_Smith_Gain}]{\includegraphics[width=0.55\textwidth]{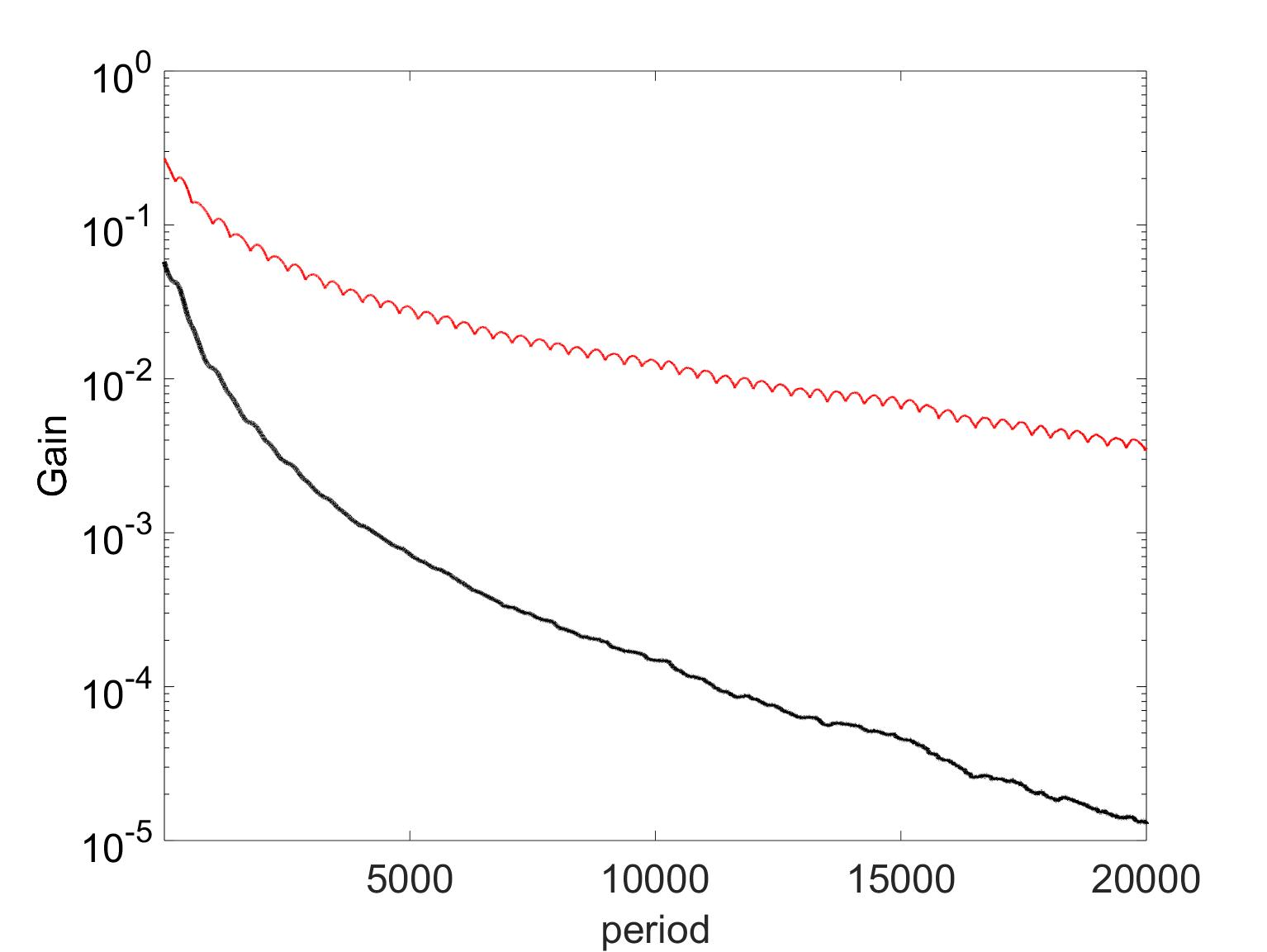}}
		
		\caption{Convergence to the unique Nash equilibrium in a good RPS and changes in the aggregate gross and net gains under a pairwise comparison dynamic (Smith). In \Cref{fig:RPS_Smith_Gain}, the aggregate net gain $G$ is shown with a black thick curve and the aggregate gross gain $\Gamma$ with a red thin curve.\label{fig:RPS_Smith}}
	\end{figure}
	Here we consider random matching in a ``good'' Rock-Paper-Scissors game \citep[Example 3.3.2]{SandholmPopText}, where a player gets payoff of $1$ from a win, $0$ from a draw and $-0.9$ from a loss. This game is strictly stable and thus the unique Nash equilibrium $\x^\ast=(x_R^\ast,x_P^\ast,x_S^\ast)=(1/3,1/3,1/3)$ is globally asymptotically stable under cost-benefit rationalizable dynamics. 
	
	In \Cref{fig:RPS_Smith_Str}, the Smith dynamic starts from $\x=(0.9,0.05,0.05)$ and indeed converges to $\x^\ast$. In \Cref{fig:RPS_Smith_Gain}, the aggregate net gain $G$ decreases over time while the aggregate gross gain $\Gamma$ has a decreasing trend but with perpetual fluctuations.
\end{exmp}

\subsection{Appendix to \Cref{sec:OrderComp}}\label{apdx:OrderComp}
Since our economically reasonability guarantees the link between static and dynamic stability while  order compatibility in \cite{Friedman91Ecta_EvolGamesEcon} does not, it would be immediate to see that order compatibility does not imply our cost-benefit rationalizablity. But, why is that so? For this, consider the example that was used in his paper to show instability under payoff monotone dynamics. Note that \cite{Friedman91Ecta_EvolGamesEcon} defines payoff monotonicity (``order compatibility" in his term) as a property of $\V$ such as i) $\dot x_a>\dot x_b\ \Leftrightarrow\ \pi_a>\pi_b$ whenever $x_a,x_b>0$ and ii) $\dot x_u=0$ whenever $x_u=0$.\footnote{This imposes monotonicity on a change in $x_a$ for each action $a$. Or, one could impose monotonicity on the growth rate of $x_a$, i.e., $\dot x_a/x_a>\dot x_b/x_b  \ \Leftrightarrow\ F_a(\x)>F_b(\x)$ (\citealt[\S4.3]{Weibull95_EvolText}; \citealt[\S 5.4.5]{SandholmPopText}).} 

\begin{exmp}[\citealt{Friedman91Ecta_EvolGamesEcon}, Counterexample 2]\label{exmp:Friedman}
	Consider a linear one-population game $\F:\Delta^3\to\R^3$ given by
	$$ \F(\x)\coloneqq \begin{pmatrix}
		-5 & -26 & 31 \\
		34	&-5 & -29\\
		-29 & 31 &-2
	\end{pmatrix}\x.$$
	This game has a unique Nash equilibrium at $\x^\ast=(1/3,1/3,1/3)$. The game is indeed strictly stable and $\x^\ast$ is a regular ESS.
	
	\begin{figure}[t!]
		\centering
		\includegraphics[scale=1]{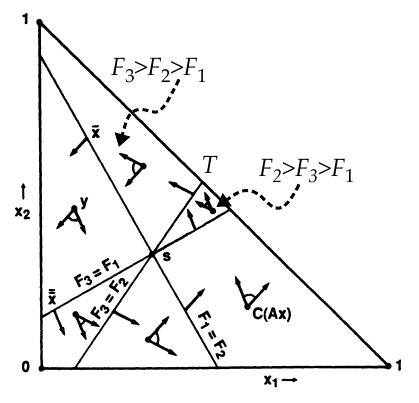}
		\caption{Phase diagram of order compatible dynamics in \Cref{exmp:Friedman}, borrowed from \citet[Figure 7]{Friedman91Ecta_EvolGamesEcon}. The labels of the payoff orderings and point $T$ (those in roman) are added to the original. The solid small arrows on the original figure indicate the possible ranges of transition vectors under order compatible dynamics.}\label{fig:FriedmanFig7}
	\end{figure}
	
	One can construct an order compatible dynamic that does not converge to $\x^\ast$ ``by always picking tangent vectors [from the cone of possible transition vectors that order compatibility permits at each point of the state] close to the outer edge of the cones [of those permissible transition vectors]'' \citep[p.655]{Friedman91Ecta_EvolGamesEcon}. Call such a dynamic an \textit{extreme diverging dynamic}. For example, consider an interior state $\x$ where $F_3(\x)> F_2(\x)> F_1(\x)$. Order compatibility requires transition vector $\dot\x$ to satisfy $\dot x_3>\dot x_2>\dot x_1$. The transition vector in an extreme diverging dynamic particularly meets $\dot x_3\approx \dot x_2>0>\dot x_1$. Thus, such a dynamic should assign to the second best action $2$ so a large net inflow as to the first best action $3$. 
	
	An extreme diverging dynamic needs to favor the \textit{second best }response action at each state almost as much as the first best, though the second best action (or the payoff ordering) may change as the state $\x$ moves. Our framework (especially, Assumption \ref{assmp:AIndp_x_gen} and independence of the availability of each action from payoff ordering of actions) allows one action to be more likely to be available than  another but this asymmetry in availability must be independent of the payoff ordering.

	In the following, we confirm that neither payoff monotonicity nor our cost-benefit rationalizablity imply each other. 

	\paragraph*{Payoff monotonicity $\nRightarrow$ cost-benefit  rationalizablity.} 
	Consider a state where  $F_3(\x)> F_2(\x)> F_1(\x)$; the transition vector in an extreme diverging dynamic must meet $\dot x_3=\dot x_2>0>\dot x_1$. Under a cost-benefit rationalizable dynamic, a transition vector $\dot\x$ should take a form such as
	\begin{align}
		\dot x_2 & = x_1 \P_{\cA 1}(\{2\})Q(F_2(\x)-F_1(\x)) &\quad - \overbrace{x_2 \P_{\cA 2}(\{\prA \mid 3\in \prA\})Q(F_3(\x)-F_2(\x)),}^{\text{gross flow from 2 to 3}} \notag\\
		\dot x_3 & = x_1 \P_{\cA 1}(\{\prA \mid 3\in \prA\}) Q(F_3(\x)-F_1(\x)) &\quad+ x_2 \P_{\cA 2}(\{\prA \mid 3\in \prA\})Q(F_3(\x)-F_2(\x)). \label{eq:Friedman_RatDyn321}
	\end{align}
	The gross flow from action 2 to action 3 (the identical second terms \textit{after }the minus/plus signs in the two equations) is non-negative. Beside, since $Q$ is non-decreasing, we have $Q(F_3(\x)-F_1(\x))\ge Q(F_2(\x)-F_1(\x))$. Therefore, $\dot x_3=\dot x_2$ holds only if  
	$$ \P_{\cA1}(\{2\}) > \P_{\cA1}(\{\prA \mid 3\in \prA\}).$$
	In other words, for an agent who has taken action 1, it must be more likely to have action 2 available but 3 not available than to have action 3 available:
	
	Now, consider an inner state where $F_2(\x)> F_3(\x)> F_1(\x)$; the transition vector $\dot\x$ in an extreme diverging dynamic should satisfy $\dot x_2=\dot x_3>0>\dot x_1$. By the same token as above, it must be the case for cost-benefit rationalizablity of $\dot x_2=\dot x_3$ that an agent who has taken action 1 is now more likely to have action $3$ available but $1$ not than having action 2 available:
	$$ \P_{\cA1}(\{3\}) > \P_{\cA1}(\{\prA \mid 2\in \prA\}).$$
	Since $\P_{\cA1}(\{b\}) \le \P_{\cA1}(\{\prA \mid b\in \prA\})$ for any  $b\in\cA\setminus\{1\}$, the above two equations contradict each other. 
	
	\paragraph*{Full rationalizability $\nRightarrow$ payoff monotonicity.} 
	Again consider a state where $F_3(\x)> F_2(\x)> F_1(\x)$. According to \eqref{eq:Friedman_RatDyn321}, a cost-benefit rationalizable dynamic may have $\dot x_2>\dot x_3$ and thus violate order compatibility if i) $\P_{\cA1}(\{2\}) \gg \P_{\cA1}(\{\prA \mid 3\in \prA\})$ and ii) $Q(F_3(\x)-F_1(\x)) \approx Q(F_2(\x)-F_1(\x)) \gg Q(F_3(\x)-F_2(\x)) \approx 0$. Below we concretely find a cost-benefit rationalizable dynamic that satisfies these two conditions. 
	
	First, we can make $\P_{\cA\cdot}$ to satisfy condition i) as well as full rationalizability (and \Cref{assmp:AExAnteavailable}): for example,
	\begin{align*}
		&\P_{\cA1}(\{2\})=3/4, &\quad& \P_{\cA1}(\{3\})=1/4; \\
		&\P_{\cA2}(\{1\})=3/4, &\quad& \P_{\cA2}(\{1,3\})=1/4; \\
		&\P_{\cA3}(\{1,2\})=3/4, &\quad& \P_{\cA3}(\{1\})=1/4.
	\end{align*}
	Note that, regardless of the current action, action 1 is available with probability 1, action 2 with 3/4 and action 3 with 1/4 as long as the action is not the current action. 
	
	To meet condition ii), focus on an inner state $\x$ near the boundary $F_3=F_2$ but a little away from $F_2=F_1$, i.e., near the intersection of  $F_3=F_2$ and $x_1+x_2=1$ (i.e., $x_3=0$); see point $T$ on \Cref{fig:FriedmanFig7}. Let function $Q(\bar q)$ be continuous (so $\cQ(\bar q)\equiv\{Q(\bar q)\}$) and strictly increasing until $\bar q$ increases to $F_3(\x)-F_2(\x)\approx 0$ but let it be constant when $\bar q$ is between $F_2(\x)-F_1(\x)$ and $F_3(\x)-F_1(\x)$. Thus, conditional switching rates from action 1 to action 3 and to action 2 are equal to each other. Yet, since action 2 is more likely to be available (alone) than action 3, there is a greater inflow from action 1 to action 2 than to action 3. 	Therefore, this dynamic fails order compatibility. But, since it is cost-benefit rationalizable, it still guarantees Nash stability in this game and thus (not only local but) global stability of the unique interior ESS.
\end{exmp}

\renewcommand\refname{References}
\begin{spacing}{1}
\bibliography{../DZbib}
\end{spacing}

\end{document}